\documentclass{conm-p-l}

\usepackage{amssymb,amsmath,amsthm,amscd,stmaryrd}
\input xy
\xyoption{all}
\usepackage[dvips]{graphicx}
\usepackage{amssymb, amsmath, amsfonts, epsfig, color, xy,amscd}

\usepackage{hyperref}

\hypersetup{
  colorlinks   = true, 
  urlcolor     = blue, 
  linkcolor    = blue, 
  citecolor   = red 
}

\usepackage[sorted,msc-links]{amsrefs}

\usepackage{amssymb}
\usepackage{amsmath,amssymb,amsbsy,amsfonts,amsthm,latexsym,amsopn,amstext, amsxtra,euscript,amscd, hyperref}
\usepackage{graphicx}
\usepackage{longtable}
\usepackage{array}

\usepackage[capitalise]{cleveref}

\crefname{thm}{Thm.}{}
\crefname{prop}{Prop.}{}
\crefname{lem}{Lem.}{}
\crefname{cor}{Cor.}{}

\usepackage{xy}
\input xy \xyoption{all}

\newtheorem{thm}{Theorem}
\newtheorem{prop}{Proposition}
\newtheorem{lem}{Lemma}
\newtheorem{cor}{Corollary}
\newtheorem{rem}{Remark}
\newtheorem{exa}{Example}
\newtheorem{alg}{Algorithm}

\newtheorem{prob}{Problem}

\theoremstyle{defi}
\newtheorem{defi}[thm]{Definition}

\def\N{\mathbb N}

\def\P{\mathbb P}
\def\F{\mathbb F}

\def\L{{\mathcal L}}
\def\O{{\mathcal O}}

\newcommand{\ch}{\mbox{\rm{char }}} 
\def\PP{\mathcal P}               
\newcommand\p{\mathfrak p}        

\DeclareMathOperator\Aut{\mathrm{Aut}}

\DeclareMathOperator\Gal{\mathrm{Gal} }

\DeclareMathOperator\mult{mult }
\DeclareMathOperator\ord{ord }
\DeclareMathOperator\dv{div }

\DeclareMathOperator\psl{\mathrm{PSL}}
\DeclareMathOperator\pgl{\mathrm{PGL}}
\DeclareMathOperator\Pic{Pic}
\DeclareMathOperator\Div{Div}

\DeclareMathOperator\diff{Diff}

\DeclareMathOperator\degrm{\mathrm{deg}}

\newcommand\X{\mathcal X}     
\newcommand\Y{\mathcal Y}     
\def\K{{\mathcal K}}          

\newcommand{\cH}{\mathcal H}

\def\a{{\alpha }}

\def\s{{\sigma }}

\def\d{\delta}

\def\t{t}
\def\r{r}

\def\iso{{\, \cong\, }}
\def\<{\langle}
\def\>{\rangle}

\def\fix{\mbox{Fix }}

\def\G{\bar G}

\def\W{\mathcal W}
\def\bC{\mathbf{C}}
\begin{document}

\title{On automorphisms of algebraic curves}

\author{A. Broughton}
\address{Department of Mathematics (Emeritus)\\
Rose-Hulman Institute of Technology\\
Terre Haute, IN 47803, USA.}
\email{brought@rose-hulman.edu}

\author{T. Shaska}
\address{Department of Mathematics and Statistics \\
Oakland University \\
Rochester, MI 48309, USA.}
\email{shaska@oakland.edu}

\author{A. Wootton}
\address{Department of Mathematics\\
University of Portland\\
Portland, OR 97203, USA.}
\email{wootton@up.edu}

\begin{abstract}
An irreducible, algebraic curve $\X_g$ of genus $g\geq 2$ defined over an algebraically closed field $k$ of characteristic $\ch k = p \geq 0$, has finite automorphism group $\Aut(\X_g)$. In this paper we describe methods of determining the list of groups $\Aut(\X_g)$ for a fixed $g\geq 2$. Moreover, equations of the corresponding families of curves are given when possible. 
\end{abstract}


\maketitle

\setcounter{tocdepth}{1}

\section{Introduction}

Algebraic curves are some of the most studied mathematical objects.  One of the first questions asked about algebraic curves was if they had any symmetries (automorphisms) and if so how many of them?  It was noticed early in XIX-century that if $g=0, 1$, then the curve could have infinitely many automorphisms.  However, if $g\geq 2$ then the set of automorphisms is finite, this was first proved by Schwartz.

Given a genus $g\geq 2$ curve $\X_g$ defined over a perfect field $k$, $\ch k = p \geq 0$, let $F:=k(\X_g)$ be its function field and denote by $\Aut(\X_g) := \Aut(F/k)$.      Historically, the main questions addressed when it comes to automorphisms of curves have been the following.

\begin{itemize}
\item[i)] What is the order  $|\Aut(\X_g)|$ ?

\item[ii)] What is the list of groups $G=\Aut(\X_g)$, for a fixed $g$?

\item[iii)] For a given $G=\Aut(\X_g)$, can one determine an equation for  $\X_g$?
\end{itemize}

The answer to each of the above questions becomes simpler if we assume that the curve $\X_g$ is a smooth, irreducible algebraic curve defined over an algebraically closed field $k$, which will be the focus of this paper.

The answer to question i) is well-known and has been long established.  When $p=0$ then $|\Aut(\X_g) | \leq 84 (g-1)$ (the so called Hurwitz bound) and when $p>0$ then $| \Aut(\X_g) | \leq  16 g^4$. In the case $p=0$ the bound is sharp and curves which achieve this bound are called Hurwitz curves.  The first example of a Hurwitz curve occurs in genus $g=3$ and is the celebrated Klein's quartic with group $\psl(2, 7)$.  The next Hurwitz curve occurs for $g=7$. It  was first discovered by Fricke and later by Macbeath \cite{Mac}.  The next two Hurwitz curves occur for $g=14, 17$.
If $p>0$ the bound is naturally higher, even though cases with $|G|>8g^3$ are very special though well-known (cf. \cref{main-thm}).

The answer to question ii) is known for $p=0$ for small genus due to recent advances in computational group theory. There is a huge amount of literature for $p=0$ from the analytic point of view via Fuchsian groups and Riemann surfaces.  In \cite{Br} Breuer computed all possible signatures of Riemann surfaces up to genus $g=48$.  The restriction on the genus was simply a restriction on the GAP's SmallGroupLibrary which included all groups of order $\leq 1000$.  This bound has been extended recently to all groups of order $< 2000$ and therefore Breuer's result can be easily extended. In \cite{kyoto} the authors, using previous results of Singerman \cite{Si}, Ries \cite{Ri},   et al.,  and their moduli dimension  gave a method of how to pick from all Breuer's signatures those signatures which occur as full automorphism groups. Moreover, all the signatures of groups of size $\geq 4 (g-1)$, also called ``large groups",  and $g\leq 10$ are given in \cite{kyoto}.  Hence, for $p=0$ we can explicitly determine the list of all groups $G=\Aut(\X_g)$ for any reasonably chosen $g\geq 2$.

In the case of $\ch k = p>0$ we are not aware of any comparable results.  However,  the list of groups for $p > 2g+1$ is the same as for $p=0$. Hence, for any given genus $g$ the fields that need to be checked are those with characteristic  $2 \leq p \leq 2g+1$.  \cref{main-thm} shows that there are very few families of curves with $|G| \geq 8g^3$. They are all of moduli dimension 0 and are superelliptic curves or cyclic curves. (A superelliptic curve (or cyclic curve) has an equation of the form $y^n=f(x)$ with some restrictions. See  \cref{sec-super} for background literature and  \cref{subsec-algextfuncfield} for a definition in terms of the automorphism group). So for all practical purposes, if we want to find a list of groups $G$ for $p>0$, then $|G| < 8g^3$ is the practical bound.

As previously remarked, the special cases of \cref{main-thm} are well known falling into the family of superelliptic curves (or cyclic curves).  Such curves are the best understood families of curves and interesting from many points of view; see \cite{bsz}, \cite{MR2367218} for more details.  In \cite{sa-1} all groups which occur as automorphism groups $\Aut(\X_g)$ over any field $k$ of $\ch k = p \neq 2$ are determined.

The well-known hyperelliptic curves are a special case of superelliptic curves.  They are the best understood among all families of algebraic curves and were the first examples used to understand Jacobians, theta-functions, etc in algebraic geometry.  In \cite{issac}, \cite{serdica}, \cite{bgg}, \cite{stich}, \cite{ShV} one can find complete details of groups that occur as automorphism groups of hyperelliptic curves, the corresponding locus in the moduli space for each such group, and invariants parameterizing some of these spaces.

The answer to question iii) is not yet completely answered. The only families of curves for which we fully understand how to write down equations and even determine such equations over a field of moduli are the superelliptic curves.  Such equations are given in \cite{s-sh} for all superelliptic curves over a field of characteristic $p\neq 2$.  In \cite{h-sh} conditions are given when such equations are defined over a field of moduli. It is still an open problem, even for characteristic $p=0$, to determine such equations for all curves.

The main goal of this paper is to give a comprehensive survey of the main results of the topic and more importantly to provide complete results for small genus $g$ or for known families of curves. We prefer the algebraic approach which makes it possible to give a more unified approach for any $p \geq 0$, but we also briefly describe the analytic approach via Riemann surfaces, Fuchsian groups, etc. We give a very brief review of definitions of the main results and refer the reader to more complete works in this topic such as \cite{Sti}, \cite{hkt}.

Throughout the paper we highlight what can be accomplished with superelliptic curves in terms of Weierstrass points, automorphism groups, and equations of curves.  We give precise lists of groups for such curves in all characteristics.

\smallskip

\noindent \textbf{Notation:}  By a \textbf{curve} $\X_g$ we always mean a smooth, irreducible, projective curve of genus $g \geq 2$ defined over a perfect field $k$.  $\Sigma_\X(k)$ will denote the set of $k$-points of $\X$. 
By $F:=k(\X_g)$ we denote its function field. For any function field $F$, a place in $F$ is denoted by $\p$ and the set of places of $F$   by $\PP_F$ (see  \cref{subsec-algextfuncfield} for definition of place).  A finite field will be denoted by $\F_q$ and the characteristic of a field $k$ by $\ch k=p$.

For a curve $\X_g$ defined over $k$, we denote by $\Aut_k(\X_g)$ the automorphism group of $\X_g$.  Since we will only consider such groups for curves defined over algebraically closed fields, then $\Aut(\X_g)$ will be used instead.
The cyclic group of order $n$ will be denoted by $C_n$, the dihedral group of order $2n$ will be $D_n$ and the symmetric and alternating groups on $n$ letters by $S_n$ and $A_n$ respectively.  In many of our results the GAP identity of a group will be used. A group $G$ with GAP identity $(m, n)$ means that $|G|=m$ and that the group is the $n$-th in the list of SmallGroup library in GAP.

\section{Algebraic curves and their function fields}

We assume that the reader is familiar with the basic definitions of field extensions.  This part is intended more to settle the notation used in the rest of the paper than as an introduction to algebraic curves. Throughout $k$ is a perfect field. For more details the reader is encouraged to see \cite{Sti} or \cite{hkt} among other places.

Let us establish some notation and basic facts about algebraic curves and their function fields.

\subsection{Algebraic curves}
The following definitions are easily extended to  any algebraic variety, but we will stick with curves. Let $k$ be a perfect field and $\X$ an algebraic curve defined over $k$.  
Then there is a homogeneous ideal  $ I_\X \subset k[X_0, X_1 ,\ldots,X_n]$ defining $\X$,  and the curve
$\X$ is irreducible if and only if $I_\X$ is a \emph{prime} ideal in  $k[X_0, X_1 ,\ldots,X_n]$.  The (homogenous) coordinate ring of $\X$ is
$ \Gamma_h (\X) := k[X_0, X_1 , \ldots , X_n]/I_\X$,
which is an integral domain. The function field of $\X$ is the quotient field of $\Gamma_h (\X)$ and denoted by $k(\X)$. Since $\X$ is an algebraic variety of dimension one, then $k(\X)$ is an algebraic function field of one variable.

Let $P=(a_0, a_1 ,\ldots , a_n)\in \X$. The ring
\[ \O_P (\X) = \{ f\in k(\X) \, | \, f \text{ is defined at } P \} \subset k(\X) \]
is a local ring with maximal ideal
\[ M_P (\X) = \{ f \in \O_P (\X) \, | \, f(P) = 0 \}   . \]
The point $P \in \X$ is a \emph{non-singular point} if the local ring $\O_P (\X)$ is a discrete valuation ring. There is a 1-1 correspondence between points $P \in \X$ and the places of $k(\X)/k$, given by  $P \mapsto M_P (\X)$.
This correspondence makes it possible to translate definitions from algebraic function fields to algebraic curves and vice-versa.


\subsection{Algebraic extensions of function fields}\label{subsec-algextfuncfield}
An algebraic function field $F/k$ of one variable over $k$ is a finite algebraic extension of $k(x)$ for some $x\in F$ which is transcendental over $k$. A \textbf{place} $\p$ of the function field $F/k$ is the maximal ideal for some valuation ring $\O$ of $F/k$.  We will denote by $\PP_F$ the set of all places of $F/k$.  Equivalently $\Sigma_\X(k)$ will denote the set of $k$-points of $\X$.

An algebraic function field $F^\prime / k^\prime$ is called an algebraic extension of $F/k$ if $F^\prime$ is an algebraic extension of $F$ and $k \subset k^\prime$.

A place $\p^\prime \in \PP_{F^\prime}$ is said to \textbf{lie over} $\p \in \PP_F$ if $\p \subset \p^\prime$. We write $\p^\prime | \p$.  In this case there exists an integer $e \geq 1$ such that $v_{\p^\prime} (x) = e \cdot v_\p (x)$, for all $x\in F$.  This integer is denoted by $e (\p^\prime  | \p) :=e$ and is called the \textbf{ramification index} of $\p^\prime$ over $\p$.
We say that $\p^\prime | \p$ is \emph{ramified} when $e (\p^\prime | \p) > 1$ and otherwise \emph{unramified}.

For any place $\p \in \PP_F$ denote by $F_\p:= \O/\p$.
The integer $f( \p^\prime | \p) := [F^\prime_{\p^\prime} : F_\p]$ is called the \emph{relative degree} of $\p^\prime | \p$.

\begin{thm}[Fundamental Equality]\label{th-fundeq}
Let $F^\prime/k^\prime$ be a finite extension of $F/k$ and $\p$ a place of $F/k$.  Let $\p_1, \dots , \p_m$ be all the places in $F^\prime/k^\prime$ lying over $\p$ and $e_i:= e(\p_i | \p)$ and $f_i :=  f(\p_i | \p)$ the relative degree of $\p_i | \p$. Then
\[ \sum_{i=1}^m e_i f_i = [F^\prime : F].\]
\end{thm}


For a place $\p \in \PP_F$ let $\O_\p^\prime$ be the integral closure of $\O_\p$ in $F^\prime$.  The complementary module over $\O_\p$ is given by $t \cdot \O_p^\prime$. Then for $\p^\prime | \p$ we define the \textbf{different exponent} of $\p^\prime$ over $\p$ as
\[ d (\p^\prime | \p) := - v_{\p^\prime} (t). \]
By \cite[Prop. 3.4.2]{Sti} $d (\p^\prime | \p) $ is well-defined and $d (\p^\prime | \p) \geq 0$.  Moreover, $d (\p^\prime | \p) =0$ for almost all $\p \in \PP_F$.  The \textbf{different divisor} is defined as
\[ \diff (F^\prime / F) := \sum_{\p \in \PP_F} \sum_{\p^\prime | \p} d(\p^\prime | \p) \cdot \p^\prime.   \]

The following well-known formula for Riemann surfaces can now be generalized for function fields as follows.

\begin{thm}[Hurwitz Genus Formula]\label{th-HGF}
Let $F/k$ be an algebraic function field of genus $g$ and $F^\prime/F$ a finite separable extension. Let $k^\prime$ denote the constant field of $F^\prime$ and $g^\prime$ the genus of $F^\prime/k^\prime$. Then,
\begin{equation}\label{e1}
2 (g^\prime-1)   =   \frac {[F^\prime : F]}    {[k^\prime: k]}   (2g-2)   +  \degrm     \diff (F^\prime /F)
\end{equation}
\end{thm}

For a proof see \cite[Thm. 3.4.13]{Sti}. A special case of the above is the following:

\begin{cor}
Let $F/k$ be a function field of genus $g$ and $x\in F\setminus k$ such that $F/k(x)$ is separable. Then,
\[ 2g-2 = -2 [F:k(x)] + \degrm \diff (F/k(x)) \]
\end{cor}

The ramification index and the different exponent are closely related, as made precise by the Dedekind theorem.

\begin{thm}[Dedekind Different Theorem] 
For all $\p^\prime|\p$ we have:

i) $d ( \p^\prime | \p) \geq e ( \p^\prime | \p) -1$.

ii) $d ( \p^\prime | \p) = e ( \p^\prime | \p) -1$ if and only if $e ( \p^\prime | \p) $ is not divisible by the $\ch k$.
\end{thm}

An extension $\p^\prime | \p$ is said to be \emph{tamely} ramified if $e (\p^\prime | \p) > 1$ and $\ch k$ does not divide $e (\p^\prime | \p)$.  If $e (\p^\prime | \p) > 1$ and $\ch k$ does  divide $e (\p^\prime | \p)$ we say that $\p^\prime | \p$ is \emph{wildly} ramified.

The extension $F^\prime / F$ is called \emph{ramified} if there is at least one place $\p \in \PP_F$ which is ramified in $F^\prime / F$. The extension  $F^\prime / F$ is called \emph{tame} if there is no place $\p \in \PP_F$ which is wildly ramified in $F^\prime/F$.

\begin{lem}
Let $F^\prime/F$ be a finite separable extension of algebraic function fields. Then

\begin{itemize}
\item[a)] $\p^\prime | \p$ is ramified if and only if $\p^\prime \leq \diff (F^\prime/F)$.  Moreover, if $\p^\prime/\p$ is ramified then:

\subitem i) $d ( \p^\prime | \p) = e ( \p^\prime | \p) -1$ if and only if $\p^\prime |\p$ is tamely ramified

\subitem ii) $d ( \p^\prime | \p) >  e ( \p^\prime | \p) -1$ if and only if $\p^\prime |\p$ is wildly ramified

\item[b)] Almost all places $\p \in \PP_F$ are unramified in $F^\prime /F$.
\end{itemize}
\end{lem}

From now on we will use the term  "curve" and its function field interchangeably, depending on the context.  It is more convenient to talk about function fields than curves in most cases.

\subsection{Divisors and the Riemann-Roch theorem}

For a given curve $\X$ defined over $k$, a divisor $D$  is called the formal finite sum 
\[D=\sum_{\p\in \Sigma_\X(k)}   z_\p \, P.\]
The set of all divisors of $\X$ is denoted by $\Div_\X (k)$.

\subsubsection{Riemann-Roch Spaces}
Define a partial ordering of elements in  $\Div_\X (k)$ as follows; $D$   is \emph{effective} ($D \geq 0$)  if $ z_\p \geq 0$ for every  $\p$,   and $D_1\geq D_2$ if $D_1-D_2\geq 0$.
The \textbf{Riemann-Roch space} associated to $D$ is
\[\L(D)=\{f\in k(\X)^*\mbox{ with } (f)\geq -D\}\cup\{0\}.\]
So the elements $x\in \L(D)$ are defined by the property that $w_\p(x)\geq -z_\p$ for all $\p\in \Sigma_\X(k)$.
Basic properties of valuations imply immediately that $\L(D)$ is a vector space over $k$. This vector space has positive dimension  if and only if there is a function $f\in k(\X)^*$ with $D+(f)\geq 0$, or equivalently,  $D\sim D_1$ with $D_1\geq 0$.

Here are some immediately obtained facts:
$\L(0)=k$ and if $\degrm(D)<0$ then $\L(D)=\{0\}$.  If $\degrm(D)=0$ then either $D$ is a principal divisor or $\L(D)=\{0\}$.
The following result is easy to prove but fundamental.
\begin{prop}
Let $D=D_1-D_2$ with $D_i\geq 0$. Then
\[\dim(\L(D))\leq \degrm (D_1)+1.\]
\end{prop}
We remark that for $D\sim D'$ we have $\L(D)\sim \L(D')$.   In particular $\L(D)$ is a finite-dimensional $k$-vector space. We follow traditional conventions and denote the dimension of $\L (D)$ by 
\begin{equation}
\ell(D):=\dim_k(\L(D)).
\end{equation}
Computing  $\ell(D)$ is a fundamental problem which is solved by the   Riemann-Roch Theorem.
A first estimate is a generalization of the proposition above:\\

For all divisors $D$ we have the inequality
\[\ell(D)\leq \degrm(D)+1.\]
For a proof one can assume that $\ell(D) > 0$ and so $D\sim D'>0$.
The important fact is that one can estimate the interval given by the inequality.

\begin{thm}[\textbf{Riemann}]  \label{riemann}
For given curve $\X$ there is a minimal number $g_\X\in \N\cup \{0\}$ such that for   all $D\in \Div_\X$ we have
\[\ell(D)\geq \degrm(D)+1-g_\X.\]
\end{thm}

\noindent For a proof see    \cite[Proposition 1.4.14]{Sti}.  Therefore,
\[ g_\X=\max \{ \degrm{D}-\ell(D)+1;\,\, D\in \Div_\X(k) \} \]
exists and is a non-negative integer    independent of $D$. The integer  $g_\X$ is called the \textbf{genus} of $\X$.

The genus does not change under constant field extensions because we have assumed that $k$ is perfect.  This can be wrong in general if the constant field of $\X$ has inseparable algebraic extensions.  There is a corollary of the theorem.

\begin{cor}
There is a number $n_\X$ such that for $\degrm(D) > n_\X$ we get equality
\[ \ell(D) =\degrm(D)+1-g_\X. \]
\end{cor}

\cref{riemann} together with its corollary is the  "Riemann part" of the Theorem of Riemann-Roch for curves. To determine $n_\X$ and to get more information about the inequality for small degrees one needs canonical divisors. \\

\subsubsection{Canonical Divisors}
Let $k(\X)$ be the function field of a curve $\X$ defined over $k$. To every $f\in k(\X)$ we attach a symbol $df$, the \emph{differential}   of $f$ lying in a $k(\X)$-vector space $\Omega(k(\X))$ generated by the symbols $df$ modulo the following relations:

For $f,g\in k(\X)$ and $\lambda\in k$ we have:

\begin{itemize}
\item[i)]  $d(\lambda  f+g)=\lambda df +dg$

\item[ii)]  $d(f\cdot g)=f dg+g df$.
\end{itemize}

\noindent The relation between derivations and differentials is given by the

\begin{defi}[Chain rule]
Let $x$ be as above and $f\in k(\X)$. Then $df = (\partial f/\partial x)  dx$.
\end{defi}

The $k(\X)$-vector space of differentials $\Omega (k(\X))$ has  dimension $1$ and it is generated by $dx$ for any $x\in k(\X)$ for which $k(\X)/k(x)$ is finite and separable. We use a well known fact from the theory of function fields $F$ in one variable.

Let $\p$ be a place of $F$, i.e. an equivalence class of discrete rank one valuations of $F$ trivial on $k$).  Then there exist a function $t_\p\in F$ with $w_\p(t_\P)=1$ and $F/k(t_\p)$ separable.

We apply this to $F=k(\X)$. For all $\p\in \Sigma_\X(k)$ we choose a function $t_\p$ as above.   For a differential $0\neq \omega\in \Omega (k(\X))$  we get
$\omega=f_\p \cdot dt_\p$.
%
The divisor $(\omega)$ is given by
\[ (\omega):= \sum_{\p\in \Sigma_\p} w_\p(f_\p)\cdot \p   \]
and  is a called a \textbf{canonical divisor} of $\X$.
The chain rule implies that this definition is independent of the choices, and the relation to differentials yields that  $(\omega)$ is a divisor.
Since $\Omega(k(\X) )$ is one-dimensional over $k(\X)$ it follows that the set of canonical divisors of $\X$  form a divisor class $\K_\X\in \Pic_\X(k)$ called the \textbf{canonical class} of $\X$.   We are now ready to formulate the Riemann-Roch Theorem.

\begin{thm}[\textbf{Riemann-Roch Theorem}]
Let $(W)$ be a canonical divisor of $\X$. For all $D\in \Div_\X(k)$ we have
$$\ell(D)=\degrm(D)+1-g_\X+\ell(W-D).$$
\end{thm}
For a proof see    \cite[Section 1.5]{Sti}.

A differential $\omega$ is \emph{holomorphic} if $(\omega)$ is an effective divisor.    The set of holomorphic differentials is a $k$-vector space denoted by $\Omega^0_\X$  which is  equal to $\L(W)$.    If we take $D=0$ respectively $D=W$ in the theorem of Riemann-Roch we get the following:

\begin{cor}
$\Omega^0_\X$ is a $g_\X$-dimensional $k$-vector space  and $\deg(W)=2g_\X-2$.
\end{cor}

For our applications  there are two further important consequences of the Riemann-Roch theorem.

\begin{cor}  The following are true:
\begin{enumerate}
\item If $\degrm(D) > 2 g_\X-2$ then $\ell(D)=\degrm(D)+1-g_\X.$

\item{In every divisor class of degree $g$ there is a positive divisor.}
\end{enumerate}
\end{cor}

\proof
Take $D$ with $\degrm(D) \geq 2g_\X -1$. Then $\degrm(W-D)\leq -1$ and therefore $\ell(W-D) =0$.  Take $D$ with $\degrm(D)=g_\X$. Then $\ell(D) =1+\ell(W-D)\geq 1$ and so there is a positive divisor in the class of $D$.
\qed


\subsection{Function fields and branched covers}

We continue to assume that $\X=\X_g$ is a smooth, irreducible, curve of genus $g=g_{\X}$ over an algebraically closed field $k$, with function field $k(\X)$.
Correspondingly, given a field $K$ of transcendence degree $1$ over $k$, then $K\backsimeq$ $k(\X^\prime)$ for some curve $\X^\prime$. As $k$ is algebraically closed, each place of
$k(\X)/k$ may be identified with a geometric point  $P$ on a smooth model of $\X$. For $P\in\X$, let $\mathcal{O}_{P}(\X)$ and  $M_{P}(\X)$ be as previously defined.
Since $\X$ is smooth at $P$ and one dimensional, there is a local parameter $z\in M_{P}(\X)$ such that $M_{P}(\X)=z\mathcal{O}_{P}(\X)$. We may write every $f\in k(\X)$
as $f=z^{e}v$, $v\in\mathcal{O}_{P}^{\ast}(\X)$. The number $e=e_{P}(f)=\nu_P(f)$ is the valuation, also called the \textbf{order of vanishing}.

We are particularly interested in the relationship between the non-constant morphisms $\pi: \X\rightarrow\Y$ of curves, which we shall call \textbf{branched coverings},
and the function fields $k(\X)$ and $k(\mathcal{Y)}$. Let $\pi:\X\rightarrow  \Y$ be a branched covering.
Then the induced map $\pi^{\ast}:k(\mathcal{Y)}\rightarrow k(\mathcal{X)}$, $f\rightarrow f\circ\pi$ is an embedding of fields, realizing $k(\mathcal{X)}$ as a finite
degree extension of $k(\mathcal{Y)}$. Conversely, if $L$ is a finite extension of $k(\mathcal{Y)}$ then there is an $\X$ and a morphism $\pi: \X\rightarrow\Y$ such that
the extension is induced via $\pi^{\ast}$. If $n = \left[ k(\mathcal{X)}:k(\mathcal{Y)}\right]$ is the degree of the extension, and the extension is separable,
then for all but finitely many points $Q\in\Y$,  $\pi^{-1}(Q)$ has $n$ points. The correspondence $\pi \leftrightarrow \pi^{\ast}$ is contravariant.

We can give a more precise statement of description of \textquotedblleft%
$\left\vert \pi^{-1}(Q)\right\vert =n$ generically\textquotedblright\ using the notion of ramification degree. Given $\pi:\X\rightarrow \Y,$ let $P\in\X,Q=\pi(P),$ and $z\in M_{P}(\X),$ $w\in$ $M_{Q}(\Y)$ be local parameters. The function $\pi^{\ast
}(w)\in$ $M_{P}(\X),$ so we can define the \textbf{ramification degree} of $\pi$ at $P$ to be $e_{\pi}(P)$ $=e_{P}(\pi^{\ast}(w)).$ The integer $e_{\pi}(P)\geq1$ and we say that $\pi$ is ramified at $P$ if $e_{\pi }(P)>1.$ Note that $e_{\pi}(P)>1$ if and only if $d\pi_{P}=0.$ We have the
following proposition.

\begin{prop}
Let $\pi:\X\rightarrow\Y$ be a branched covering of degree $n$. Then for $Q\in\Y$ we have
\begin{equation}
n= {\displaystyle\sum\limits_{\pi(P)=Q}}  e_{\pi}(P) \label{eq-fibrecount}%
\end{equation}
If $\pi: \X\rightarrow\Y$ is induced by a $G$-action and the field extension is separable, then for each $Q\in\Y$ number of points lying over $Q$ is $\left\vert G\right\vert /\left\vert G_{P}\right\vert $ for
any $P\in\pi^{-1}(Q).$ It follows then that
\begin{equation}
n    =e_{\pi}(P)\frac{\left\vert G\right\vert }{\left\vert G_{P}\right\vert },\label{eq-Gfibrecount}
\end{equation}
where  $e_{\pi}(P)    =\left\vert G_{P}\right\vert$
(See equation \eqref{eq-stab}).
\end{prop}

Equation \eqref{eq-fibrecount} is simply \cref{th-fundeq}.

\begin{exa}\label{ex-insep}
Let $k$ $=\overline{\F}_{p}$, and consider the Frobenius morphism $\pi: \P^{1}(k)\rightarrow\P^{1}(k),$ given by $x\rightarrow x^{p}.$ The map is an injective but not invertible morphism. The induced maps on fields $k(y)\rightarrow k(x)$ is given by $y=x^{p},$ and has degree $p$. By direct calculation it can be shown that $e_{\pi}(x)=p$ for all $x\in\P^{1}(k),$ or alternatively $d\pi=0$. Thus equation \eqref{eq-fibrecount} is satisfied. This strange behavior is linked to the fact that the extension $k(y)=k(x^{p})\rightarrow k(x)$ is purely inseparable. 
\end{exa}

From now on we will consider only separable extensions.

A branched cover $\pi:\X\rightarrow\Y$ is called \textbf{tamely ramified} if all branching orders $e_{\pi}(P)$,  $P\in \X$ are relatively prime to the characteristic $p$. In characteristic $p=0$ all branched covers are considered to be tamely ramified. A cover is \textbf{wildly ramified} if it is not tamely ramified.
 
%

\subsection{Automorphism groups, $G$-actions, stabilizers}


Let $\X_g$ be an irreducible and non-singular algebraic curve defined over $k$ of genus $g\geq 2$. We denote its function field by $F:=k(\X_g)$.  The automorphism group of $\X_g$ is the group $G:=\Aut(F/k)$ (i.e., all field automorphisms of $F$ fixing $k$). We will denote it by $\Aut_k (\X_g)$.  When $k$ is algebraically closed then we will simply use $\Aut (\X_g)$.  
The rest of this paper will focus on determining $\Aut (\X_g)$, for any given $g\geq 2$.

We say that  a finite group $G$ acts (birationally, conformally) on $\X$ if there is a monomorphism $\epsilon : G \rightarrow \Aut(\X)$. There is an induced action on function fields given by
$\epsilon^{\ast} : h\rightarrow \left(  h^{-1}\right)  ^{\ast}$.
Assuming $G$ is finite, then the field of invariant functions $k(\X)^{G}$ is a subfield such that $k(\X)$ is an extension of $k(\X)^{G}$ of degree $\left\vert G\right\vert$. The subfield $k(\X)^{G}$ corresponds to some $k(\Y)$ and there is a morphism $\pi_{G} : \X \rightarrow \Y.$ We denote $\Y$ by $\X/G.$ It can be shown that $G$ acts transitively on the fibers of $\pi_{G}$, so that $\Y$ is an orbit space, as a set, so that degree of $\pi_{G}$ is $\left\vert G\right\vert$. The set of \textbf{branch points} or \textbf{branch locus} is denoted $\mathcal{B}_{G}.$ The covering $\pi_{G}: \X^{\circ} \rightarrow \Y^{\circ}$, where $\Y^{\circ}=\Y-\mathcal{B}_{G}$ and $\X^{\circ}=\X-\pi_{G}^{-1}(\mathcal{B}_{G})$ is an unramified Galois covering of affine curves.

A first difference that we see in positive characteristic is the structure and action of stabilizers. Given $G$ acting on $\X$ and $P\in\X,$ the \textbf{stabilizer} at $P$ is defined by
\begin{equation}
 G_{P} =\left\{  g\in G : gP=P\right\},  \label{eq-stab}
\end{equation}
also known as the decomposition group. In the characteristic $0$ case $G_{P}$ is cyclic, and acts faithfully on the tangent space $T_{P}(\X).$ This may fail in the $p>0$ case.

The action of stabilizers on $M_{P}(\X)/ M_{P}^{s+1}(\X)$ is faithful for some $s$. Indeed, if $x$ is a local parameter at $P$, and $g^{\ast}$ is the identity on $M_{P}(\X)/M_{P}^{s+1}(\X)$ then $g^{\ast} x - x\in M_{P}^{s+1}(\X)$ for all $s$. It follows that $g^{\ast}x-x=0$ as ${\displaystyle\bigcap\limits_{s\geq0}} M_{P}^{s+1}(\X)=\left\{  0\right\}$.
Now let $f\in$ $\mathcal{O}_{P}(\X)$ be arbitrary. By considering a Taylor series expansion in $x$ at $P$ we see that $g^{\ast}f=f$. It follows that $g^{\ast}$ is  trivial
on all of $k(\X)$ and hence $g$ is trivial on all of $\X$. From the format of the matrix for $g^{\ast}$ it follows that the action on
$GL_k\left(M_{P}(\X)/M_{P}^{2}(\X)\right) \backsimeq k^{\ast}$ is trivial on elements of order $p$, so that the map $G_{P}\rightarrow k^{\ast}$ has cyclic image and the kernel
is a $p$-group.
It then follows that $G_{P} \backsimeq  C_{m}\ltimes H$ where $H$ is a $p$-group and that $e_{\pi}(P)$ $=m\times q$ where $q$ is some $p$-power. For $g\in$ $G_{P}$ we call the image of $g^{\ast}$ in $k^{\ast}$ the \textbf{rotation number} of $g$ at $P$.

\subsection{Cyclic n-gonal curves}

A cyclic $n$-gonal curve has a cyclic $G=C_n$ action for which the genus of $\X/G$ is zero. A defining equation is of the form $y^{n}=f(x),$ where $f(x)$ is a rational function. The curve can be put in a
canonical form:
\[
y^{n}=(x-a_{1})^{n_{1}}  \cdots     (x-a_{t})^{n_{t}},
\]
where $n,n_{1},\ldots,n_{t}\in\mathbb{Z}^{+},$ $a_{1},\ldots,a_{t}\in k$ satisfy:

\begin{enumerate}
\item $a_{1},\ldots,a_{t}$ are distinct,

\item $0<n_{i}<n,$

\item $n$ divides $n_{1}+\cdots+n_{t}$

\item $\gcd(n_{1},\ldots,n_{t})=1.$
\end{enumerate}

\noindent The conditions 1 and 2 simplify the model and eliminate degeneracies, condition 3 ensures that curve is not ramified over $\infty,$ and condition 4 ensures that curve is irreducible. The curve needs to be normalized to make it smooth. There is a $G$-action, where $G = U_{n}=\left\{  u\in k:u^{n}=1\right\}$, given by $(x,y)\rightarrow(x,uy).$ We must assume the $p$ does not divide $n$ otherwise the field extension will be inseparable. The quotient map is $(x,y)\rightarrow x,$ and is branched over $a_{1},\ldots,a_{t}.$ The local equation is $y^{n}=b(x)(x-a_{i})^{n_{i}},$ where $b(a_{i})\neq 0.$ There are $d_{i}$ $=\gcd(n,n_{i})$ branches of the curve meeting at $(a_{i},0),$ and in the normalization the ramification degree is $m_{i}=n/d_{i}.$ For a discussion of the $p=0$ case, which extends in many ways to the $p>0$ case see \cite{B1},\cite{B2}. In the complex case, the vector $(n_1, \ldots, n_t)$ is a generating vector (see \cref{subsec-FGS}) for $G=C_{n}$ by the conditions $2,3,4.$ In the $p>0$ case $(n_1, \ldots, n_t)$ is a classifier for the $C_{n}$ action and may be thought of as a generating vector for Galois covers of $\Y=\P^{1}(k),$ ramified exactly over $\mathcal{B}=\{a_{1},\ldots,a_{t}\},$ with Galois group $G=C_{n}$. It plays the role of generating vector, though it is not constructed topologically. We will look at such curves in more detail in \cref{sec-super}.

\section{Weierstrass Gap Theorem and Weierstrass points}\label{subsec-weierstrass}
We assume that the reader is familiar with basic definitions on divisors on curves. For a short introduction see \cite{Frey} in this volume or \cite{Sti}.

Let $P$ be a point on $\X_g$ and consider the vector spaces $\L(nP)$ for $n=0,1,\dots,2g-1$.  These vector spaces contain functions with poles only at $P$ up to a specific order.  This leads to a chain of inclusions \[ \L(0)\subseteq \L(P) \subseteq \L(2P) \subseteq \dots \subseteq \L((2g-1)P)\] with a corresponding non-decreasing sequence of dimensions \[ \ell(0) \leq \ell(P) \leq \ell(2P) \leq \dots \leq \ell((2g-1)P).\]
The following proposition shows that the dimension goes up by at most 1 in each step.

\begin{prop}\label{prop:dimension-increase-by-1} For any $n>0$,
\[ \ell((n-1)P)\leq \ell(nP)\leq \ell((n-1)P)+1.\]
\end{prop}

\proof
It suffices to show $\ell(nP)\leq \ell((n-1)P)+1$.  To do this, suppose $f_1, f_2\in\L(nP)\setminus\L((n-1)P).$  Since $f_1$ and $f_2$ have the same pole order at $P$, using the series expansions of $f_1$ and $f_2$ with a local coordinate, one can find a linear combination of $f_1$ and $f_2$ to eliminate their leading terms.  That is, there are constants $c_1, c_2\in k$ such that $c_1f_1+c_2f_2$ has a strictly smaller pole order at $P$, so $c_1f_1+c_2f_2\in\L((n-1)P)$.  Then $f_2$ is in the vector space generated by a basis of $\L((n-1)P)$ along with $f_1$.  Since this is true for any two functions $f_1,f_2$, we conclude $\ell(nP)\leq \ell((n-1)P)+1$, as desired.\qed

For any integer $n>0$, we call $n$ a \textbf{Weierstrass gap number of $P$} if $\ell(nP)=\ell((n-1)P)$; that is, if there is no function $f\in k(\X_g)^\times$ such that $(f)_\infty=nP$.  Weierstrass stated and proved the Gap theorem, or \textbf{L\"uckensatz}, on gap numbers in the 19th century, likely in the 1860s.

\begin{thm}[The Weierstrass Gap Theorem]
For any point $P$, there are exactly $g$ gap numbers $\alpha_i(P)$ with \[1=\alpha_1(P)<\alpha_2(P)<\cdots<\alpha_g(P)\leq2g-1.\]
\end{thm}
This theorem is a special case of the Noether Gap theorem, which we state and prove below.

The set of gap numbers, denoted by $G_P$,  forms the \textbf{Weierstrass gap sequence} for $P$.  The non-gap numbers form a semi-group under addition since they correspond to pole orders of functions. For any curve the gap sequence is the same for all points  with finitely many exceptional points. All curves in characteristic 0 and most curves in positive characteristic have the \textbf{classical gap sequence} $\{1,2,\dots,g\}$ for the generic gap sequence. In this section we discuss only curves with a classical gap sequence. See \cite{Sch} for examples of curves with a non-classical gap sequence.

\begin{defi}[Weierstrass point]
If the gap sequence at $P$ is anything other than $\{1,2,\dots,g\}$, then $P$ is called a \textbf{Weierstrass point}.
\end{defi}
Equivalently, $P$ is a Weierstrass point if $\ell(gP)>1$; that is, if there is a function $f$ with $(f)_\infty=mP$ for some $m$ with $1<m\leq g$.
The following was proved by F. K. Schmidt (1939).

\begin{thm}[Schmidt \cite{Sch}]
Every algebraic curve of genus $g\geq 2$ has at least one Weierstrass point.
\end{thm}

The notion of gaps can be generalized, which we briefly describe.  Let $P_1,P_2,\dots,$ be a sequence of (not necessarily distinct) points on $\X_g$.  Let $D_0=0$ and, for $n\geq1$, let $D_n=D_{n-1}+P_n$.  One constructs a similar sequence of vector spaces 
\[\L(D_0)\subseteq\L(D_1)\subseteq\L(D_2)\subseteq\dots\subseteq\L(D_n)\subseteq\cdots\] 
with a corresponding non-decreasing sequence of dimensions 
\[\ell(D_0)\le\ell(D_1)\le\ell(D_2)\le\dots\le\ell(D_n)\le\cdots.\]  
If $\ell(D_n)=\ell(D_{n-1})$, then $n$ is a \textbf{Noether gap number} of the sequence $P_1, P_2, \dots .$

\begin{thm}[The Noether GAP Theorem]
For any sequence $P_1,P_2, \dots$, there are exactly $g$ Noether gap numbers $n_i$ with \[1=n_1<n_2<\dots<n_g\leq 2g-1.\]
\end{thm}

\proof In analog with \cref{prop:dimension-increase-by-1}, one can show the dimension goes up by at most 1 in each step; that is, \[ \ell(D_{n-1})\leq \ell(D_n)\leq \ell(D_{n-1})+1\] for all $n>0$.  First, note that the Riemann-Roch theorem is an equality for $n>2g-1$, so the dimension goes up by 1 in each step, so there are no gap numbers greater than $2g-1$.

Now, consider the chain $\L(D_0)\subseteq\dots\subseteq\L(D_{2g-1})$.  By Riemann-Roch, $\ell(D_0)=1$ and $\ell(D_{2g-1})=g$, so in this chain of vector spaces, the dimension must increase by 1 exactly $g-1$ times in $2g-1$ steps.  Thus, for $n\in\{1,2,\dots,2g-1\}$,  there are $g$ values of $n$ such that $\ell(D_n)=\ell(D_{n-1})$.  These $g$ values are the Noether gap numbers.
\qed


For a complete treatment of Weierstrass points and their weights for $p=0$ see \cite{w-2}.

\subsection{Weierstrass points via holomorphic differentials}
Continuing with a point $P$ on a curve $\X_g$, recall that $n$ is a gap number precisely when $\ell(nP)=\ell((n-1)P)$.  By Riemann-Roch, this occurs exactly when \[\ell(K-(n-1)P)-\ell(K-nP)=1\] for a canonical divisor $K$, which is the divisor associated to some differential $dx$.  Thus there is $f\in k(\X_g)^\times$ such that $(f)+K-(n-1)P\geq0$ and $(f)+K-nP\not\geq0$, which implies that $\ord_P(f\cdot dx)=n-1$.  Since $(f)+K\geq(n-1)P\geq0$ (for $n\geq1$), $n$ is a gap number of $P$ exactly when there is a holomorphic differential $f\cdot dx$ such that $\ord_P(f\cdot dx)=n-1$.

For $H^0(\X_g,\Omega^1)$ the space of holomorphic differentials on $\X_g$, by Riemann-Roch, the dimension of $H^0(\X_g,\Omega^1)$ is $g$.  Let $\{\psi_i\}$, for $i=1,\dots, g$, be a basis, chosen in such a way that
\[\ord_P(\psi_1)<\ord_P(\psi_2)<\cdots<\ord_P(\psi_g).\]
Let $n_i=\ord_P(\psi_i)+1$.
The \textbf{1-gap sequence at $P$} is $\{n_1,n_2,\dots,n_g\}$.
We then have the following equivalent definition of a Weierstrass point.

\begin{defi}[Weierstrass point]
If the 1-gap sequence at $P$ is anything other than $\{1,2,\dots,g\}$, then $P$ is a Weierstrass point.
\end{defi}

With this formulation, we see $P$ is a Weierstrass point exactly when there is a holomorphic differential $f\cdot dx$ with $\ord_P(f\cdot dx)\geq g$.

\begin{defi}[Weierstrass weight]
The \textbf{Weierstrass weight} of a point $P$ is \[w(P)=\sum_{i=1}^g (n_i-i).\]
\end{defi}

In particular, $P$ is a Weierstrass point if and only if $w(P)>0$.

\subsection{Bounds for weights of Weierstrass points}
Suppose $\X_g$ is a curve of genus $g\geq1$, $P\in\X_g$, and consider the 1-gap sequence of $P$ $\{n_1,n_2,\dots,n_g\}$.  We will refer to the non-gap sequence of $P$ as the complement of this set within the set $\{1,2,\dots,2g\}$.  That is, the non-gap sequence is the sequence $\{\alpha_1,\dots,\alpha_g\}$ where \[1<\alpha_1<\dots<\alpha_g=2g.\]
\begin{prop}\label{prop:sums-of-non-gaps}
For each integer $j$ with $0<j<g$, $\alpha_j+\alpha_{g-j}\geq 2g$.
\end{prop}

\begin{prop}\label{prop:max-w-weight}
For $P\in\X_g$, $w(P)\leq g(g-1)/2$, with equality if and only if $P$ is a branch point on a hyperelliptic curve $\X_g$.
\end{prop}

\begin{cor}\label{cor:number-w-pts}
For a curve of genus $g\geq2$, there are between $2g+2$ and $g^3-g$ Weierstrass points.  The lower bound of $2g+2$ occurs only in the hyperelliptic case.
\end{cor}
\proof The total weight of the Weierstrass points is $g^3-g$.  In \cref{prop:max-w-weight}, we see that the maximum weight of a point is $g(g-1)/2$, which occurs in the hyperelliptic case.  Thus, there must be at least $\dfrac{g^3-g}{g(g-1)/2} = 2g+2$ Weierstrass points.  On the other hand, the minimum weight of a point is 1, so there are at most $g^3-g$ Weierstrass points.
\qed

\begin{thm}[Weierstrass Normal Form] 
Let $\X$ be an irreducible curve of genus $g\geq 1$ defined over $k$.  For a place $\p \in k (\X)$ let $m$ be the first non-gap at $\p$ and $n$ be the least non-gap which is prime to $m$.  Then

i) $\X$ has affine equation
\begin{equation}\label{W-form}
 f(x, y) = y^m+ u_1 (x) y^{m-1} + \dots + u_{m-1} (x) y + u_m (x),
\end{equation}
where $u_i (x) \in k[x]$, $\degrm u_m =n$ and $\degrm u_i (x) < in/m$, for $i=1, \dots , m-1$.

ii) $p=(x, y)$ is a generic point of $\X$ for $\dv (x)_\infty = m\p$ and $\dv (y)_\infty = n\p$.

iii) the branch of $\X$ associated to $\p$ is the unique branch of $\X$ with centre at $y_\infty$.
\end{thm}

The following is a well-known result.

\begin{lem}
Let $\X$ be a genus $g$ curve in Weierstrass normal form as in \cref{W-form}. Then

i) $g \leq \frac 1 2 (n-1) (m-1)$

ii) If no point of $\X$ other than $y_\infty$ is singular, then $g=\frac 1 2 (n-1) (m-1)$
\end{lem}

For more details about Weierstrass points and their weights see \cite{w-1}, \cite{w-2} where Weierstrass points of superelliptic curves are studied.

\section{Automorphisms of curves}

Let $\X$ be an irreducible and non-singular algebraic curve defined over $k$. We denote its function field by $F:=k(\X)$.  The automorphism group of $\X$ is the group $G:=\Aut(F/k)$ (i.e., all field automorphisms of $F$ fixing $k$).

\subsection{The action of $k$-automorphisms on places}
$G$ acts on the places of $F/k$.  Since there is a 1-1 correspondence between places of $F/k$ and points of $\X$, this action naturally extends to the points of $\X$. For $\alpha \in G$ and $P \in \X$, we denote its image under $\alpha$ by $P^\alpha$.
In a natural way we extend this $G$-action to $\Div_k (\X)$. Let $D \in \Div_k (\X)$, say $D= \sum n_P \cdot P$. Then
\[ D^\alpha = \sum n_p \cdot P^\alpha.\]

\begin{lem}\label{lem-1} 
$G$ acts on the set $\W$ of Weierstrass points.
\end{lem}

\proof
The set $\W$ of  Weierstrass points do not depend on the choice of the local coordinate and so it is invariant under any $\sigma \in \Aut(\X_g)$.
\qed

Hence, in order to determine the automorphism group we can just study the action of the group on the set of Weierstrass point of the curve.
Then we have the following.

\begin{prop} 
Let $\a \in \Aut(\X)$ be a non-identity element. Then $\a$ has at most $2g+2$ fixed places.
\end{prop}

\proof
Let $\a$ be a non-trivial element of  $\Aut(F/k)$.  Since $\a$ is not the identity, there is some place $\p \in \PP_F$ not fixed by $\a$.
Take $g+1$ distinct places $\p_1, \dots , \p_{g+1} $ in $\PP_F$ such that $D = \p_1 + \cdots +\p_{g+1}$ and $D^\a$ share no place. By \cite[Thm. 6.82]{hkt}  there is $z\in F\setminus k$ such that $\dv (z)_\infty = D$.  Then consider $w= z - \a (z)$.  Since $z$ and $\alpha (z)$ have different poles then $w\neq 0$.  Hence, $w$ has exactly $2g+2$ poles. Then $w$ has exactly $2g+2$ zeroes. But every fixed place of $\a$ is a zero of $w$.  Hence $\a$ has at most $2g+2$ fixed places.
\qed

Let $\W$ be the set of Weierstrass points.  From \cref{cor:number-w-pts} we know that $\W$ is finite. Since for every $\a  \in \Aut(\X)$, from \cref{lem-1} we have $\a(\W)=\W$.  Then we have the following; see \cite[Thm. 11.24]{hkt} for the proof.

\begin{thm}  
Let $\X$ be a genus $g\geq 2$ irreducible, non-hyperelliptic curve defined over $k$ such that $\ch k = p$ and $\a \in \Aut(\X)$. If $p=0$ or $p > 2g-2$ then $\a$ has finite order.
\end{thm}

Then we have the following;  see \cite[Lemma 11.25]{hkt}.

\begin{lem}   
If $p=0$ and $g\geq 2$ then every automorphism is finite.
\end{lem}

In the case of $p=0$ Hurwitz \cite{Hu} showed $|\a| \leq 10 (g-1)$. In 1895, Wiman improved this bound to be $ |\a| \le 2(2 g+1)$ and showed this is best possible. If $|\a|$ is a prime then $|\a| \leq 2g +1$. Homma \cite{Ho} (1980) shows that this bound is achieved for a prime $q \neq p$  if and only if the curve is birationally
equivalent to
\[ y^{m-s} (y-1)^s = x^q, \quad for \quad 1 \leq s < m \leq g +1. \]

If $p > 0$ then we have the following; see \cite[Thm.~11.34]{hkt}.

\begin{thm}  
Let $\X$ be a genus $g\geq 2$, irreducible curve defined over $k$, with $\ch k =p >0$ and $\alpha \in \Aut(\X)$ which fixes a place $\p \in \PP_F$.  Then the order of $\alpha $ is bounded by
\[ | \alpha | \leq 2 p (g+1) (2g+1)^2. \]
\end{thm}

\subsection{Finiteness of $\Aut(\X)$}

The main difference for $g=0, 1$ and $g\geq 2$ is that for $g\geq 2$ the automorphism group is a finite group.  This result was proved first by Schmid (1938).

\begin{thm}[Schmid \cite{Schmid}]
Let $\X$ be a genus $g\geq 2$, irreducible curve defined over a field $k$, $\ch k = p \geq 0$. Then $\Aut(\X)$ is finite.
\end{thm}

\subsubsection{Characteristic $p=0$}
As an immediate consequence of the Hurwitz theorem for $\ch k=p=0$ we have that
\[ |\Aut(\X)| \leq 84 (g-1).\]
Curves which obtain this bound are called \emph{Hurwitz curves}.  Klein's quartic is the only Hurwitz curve of genus $g\le3$. Fricke showed that the next Hurwitz group occurs for $g=7$ and has order 504. Its group is $SL(2,8)$, and an equation for it was computed by Macbeath \cite{Mac} in 1965. Further Hurwitz curves occur for $g=14$ and $g=17$ (and for no other values of $g\le19$).

For a fixed $g\geq 2$ denote by $N(g)$ the maximum of the $|\Aut(\X_g)|$. Accola \cite{Ac1} and Maclachlan
\cite{Mc1} independently show that $N(g) \geq 8 (g+1)$ and this bound is sharp for infinitely many $g$'s. If
$g$ is divisible by 3 then $N(g) \geq 8 (g+3)$.

The following terminology is rather standard. We say $G\le\Aut(\X_g)$ is a {\bf large automorphism group} in
genus $g$ if
$$|G|\ \ >\ \ 4 (g-1).$$
Then the quotient of $\X_g$ by $G$ is a curve of genus $0$, and the number of points of this quotient ramified
in $\X_g$ is 3 or 4 (see \cite{Br}, Lemma 3.18, or \cite{FK}, pages 258-260). Singerman \cite{S1} (1974) shows
that Riemann surfaces with large cyclic, Abelian, or Hurwitz groups are symmetric (admit an anti-conformal involution).
Kulkarni \cite{Ku}(1997) classifies Riemann surfaces admitting large cyclic automorphism groups and works out
full automorphism groups of these surfaces. Matsuno \cite{Mt}(1999) investigates the Galois covering of the
projective line from compact Riemann surfaces with large automorphism groups.

\subsubsection{Characteristic $p>0$}

In the case of positive characteristic the bound is higher due to possibly wild ramifications.  The following was proved by Stichtenoth (1973) by extending previous results of P. Roquette and others.

\begin{thm}[Stichtenoth \cite{Stichtenoth}]
Let $\X$ be a genus $g\geq 2$, irreducible curve defined over a field $k$, $\ch k = p>0$. Then
\[ |\Aut(\X)| < 16 \cdot g^4, \]
unless $\X$ is the curve with equation
\[ y^{p^n} + y= x^{p^{n+1}}, \]
in which case it has genus $g= \frac 1 2 p^n (p^n-1) $ and $ |\Aut(\X)| = p^{3n} (p^{3n} +1) (p^{2n}-1)$.
\end{thm}

Hence we have a bound for curves of genus $g\geq2$ even in characteristic $p>0$.  It turns out that all curves with large groups of automorphisms are special curves.  So getting ``better" bounds while such curves are left out has always been interesting.  There is a huge amount of literature on this topic due to the interest of such bounds in coding theory. Perhaps the following theorem, which  is due to Henn \cite{Henn}, provides a better bound if the following four families of curves are left out.
As Henn points out in a footnote,   this result may be sharpened to show that the order of $\Aut(\X)$  is less than $3\cdot (2g)^{5/2}$ except when $k (\X)$  belongs to one of five types of function fields. 
   
\begin{thm}[Henn \cite{Henn}] \label{main-thm}  
Let $\X$ be an irreducible curve of  genus $g\geq 2$.  If $|G| \geq 8g^3$, then $\X$ is isomorphic to one of the following:

i) The hyperelliptic curve
\[ y^2+y + x^{2^k+1} =0, \]
defined over a field of characteristic $p=2$.  In this case the genus is $g=2^{k-1}$ and $|G|=2^{2k+1} (2^k+1)$.

ii)  The hyperelliptic curve
\[  y^2=x^q- x,\]
defined over a field of characteristic $p>2$ such that $q$ is a power of $p$.   In this case $g=\frac 1 2 (q-1)$ and the reduced group $\bar G$ is isomorphic to $\psl(2, q)$ or $\pgl(2, q)$.

iii) The Hermitian curve
\[ y^q + y = x^{q+1}, \]
defined over a field of characteristic $p\geq 2$ such that $q$ is a power of $p$. In this case $g= \frac 1 2 (q^2-q)$ and $G$ is isomorphic to $PSU(3, q)$ or $PGU(3, q)$.

iv) The curve
\[ y^q+y = x^{q_0} (x^q+x), \]
for $p=2$, $q_0=2^r$, and $q=2q_0^2$.  In this case, $g= q_0 (q-1)$ and $G \iso Sz (q)$.
\end{thm}

\subsubsection{Bounds for $p=0$}

Next we will consider some special curves, namely  superelliptic or cyclic  curves.    Most of the groups of automorphisms for a fixed $g\geq 2$ correspond to such curves and they have large automorphism groups.  A more general interesting problem which we do not consider here is the following.

\begin{prob}
Let $\X$ be an irreducible, smooth algebraic curve of genus $g\geq 2$ which is not superelliptic.  Determine a bound for $|\Aut(\X) |$.
\end{prob}

\section{Superelliptic curves}\label{sec-super}
There is a lot of literature published on superelliptic curves in the last two decades. We mostly follow the terminology from \cite{super-1}, \cite{super-2}, or \cite{h-q-sh}.

Let $k$ be a field of characteristic $p$ and $\X_g$ be a genus $g\geq 2$ cyclic curve given by the equation $y^n=f(x)$ for some $f \in k[x]$.  Then $k(x)$ is degree $n$ genus zero subfield of $K=k(\X)$. Let $G=\Aut(K/k)$. We are going to assume that $C_n:=\Gal(K/k(x))=\langle w \rangle$, with $w^n=1$ and which satisfies $\langle w \rangle \lhd G$. (This is assumed in the literature cited above, and is the definition we take in \cref{subsec-revisit}.) Then, the reduced group $\G:=G/C_n$ acts on $\P^1(k)$, and so satisfies $\G \leq \pgl(2,k)$. Hence $\G$ is isomorphic to one of the following: $C_m $, $ D_m $, $ A_4 $, $ S_4 $, $A_5$, \emph{semidirect product of an elementary  Abelian  group  with  cyclic  group}, $\psl(2,q)$ and $\pgl(2,q)$, see \cites{VM}.    

The group $\G$ acts on $k(x)$ naturally. The fixed field is a genus 0 field, say $k(z)$. Thus $z$ is a degree $|\G|$ rational function in $x$, say $z=\phi(x)$, yielding a Galois branched covering $\phi: \P^1 \to \P^1$.

Let $\phi_0:\X_g \to \P^1$ be the cover which corresponds to the degree $n$ extension $K/k(x)$. Then $\Phi := \phi \circ \phi_0$ has monodromy group $G:=\Aut(\X_g)$. From basic covering theory, the group $G$ is embedded in the group $S_l$ where $l= \degrm \Phi$. There is an $r$-tuple $\overline{\sigma}:= (\sigma_1, \dots ,\sigma_r)$, where $\sigma_i \in S_l$ such that $\sigma_1, \dots ,\sigma_r$ generate $G$ and $\sigma_1 \cdots \sigma_r=1$. The signature of $\Phi$ is an $r$-tuple of conjugacy classes $\bC:=(C_1, \dots ,C_r)$ in $S_l$ such that $C_i$ is the conjugacy class of $\sigma_i$. We use the notation $n$ to denote the conjugacy class of permutations which is cycle of length $n$. Using the signature of $\phi: \P^1 \to \P^1$ one finds out the signature of $\Phi: \X_g \to \P^1$ for any given $g\geq 2 $ and $G$.

Let $E$ be the fixed field of G, the Hurwitz genus formula states that
\begin{equation}\label{e3}
2(g_K-1)=2(g_E-1)|G|+\degrm(\diff_{K/E})
\end{equation}
with $g_K$ and $g_E$ the genera of $K$ and $E$ respectively and $\diff_{K/E}$ the different of $K/E$. Let $\overline{\p}_1, \overline{\p}_2, ..., \overline{\p}_r$ be ramified primes of $E$. If we set $d_i= \degrm (\overline{\p}_i)$,  let $e_i$ be the ramification index of the $\overline{\p}_i$ and  $\beta_i$ be the exponent of $\overline{\p}_i$ in $\diff_{K/E}$. Hence, the above equation  may be written as
\begin{equation}\label{e2}
2(g_K-1)=2(g_E-1)|G|+|G| \sum_{i=1}^{r} \frac{\beta_i}{e_i} d_i
\end{equation}

If $\overline{\p}_i$ is tamely ramified then $\beta_i=e_i-1$ or if $\overline{\p}_i$ is wildly ramified then $\beta_i=e_i^*q_i+q_i-2$ with $e_i=e_i^*q_i$, $e_i^*$ relatively prime to $p$, $q_i$ a power of $p$ and $e_i^*|q_i-1$. For fixed $G$, $\bC$ the family of covers $\Phi: \X_g \to \P^1$ is a Hurwitz space $\cH (G,\bC)$. The space  $\cH(G,\bC)$ is an irreducible algebraic variety of dimension $\d (G,\bC)$. Using equation~\eqref{e2} and signature $\bC$ one can find out the dimension for each $G$.

Next we want to determine the cover $z=\phi(x): \P^1 \to \P^1$ for all characteristics. Notice that the case of $\ch k =0$ is already worked out in  \cite{serdica} and for $\ch k = p \neq 2$ is done in \cite{sa-1} which we are following here.

We define a semidirect product of elementary Abelian group  with  cyclic group as follows.
\[
K_m:=\left\langle\left\{ \sigma_a, \t |   \,   a \in U_m   \right\}\right\rangle
\]
where $\t(x)=\xi^2x$, $\sigma_a(x)=x+a,$ for each $a \in U_m$,
\[
U_m :=\{a \in k \, | \,  (a\prod_{j=0}^{\frac{p^t-1}{m}-1}(a^m-b_j))=0\}
\]
$b_j \in \F^*_q$, $m|p^t-1$ and $\xi$ is a primitive $2m$-th root of unity. Obviously $U_m$ is a subgroup of the additive group of $k$.

\begin{lem}\label{l1}
Let $k$ be an algebraically closed field of characteristic $p$, $\G$ be a finite subgroup of $\pgl(2,k)$ acting on the field $k(x)$. Then, $\G$ is isomorphic to one of the following groups $C_m$, $D_m$, $A_4$, $S_4$, $A_5$, $U= C_p^t$, $K_m$, $\psl(2,q)$ and $\pgl(2,q)$, where $q=p^f$ and $(m,p)=1$. Moreover, the fixed subfield $k(x)^{\G}=k(z)$ is given by \cref{t1}, where $\alpha =\frac{q(q-1)}{2}, \quad \beta= \frac{q+1}{2}$. The subgroup  $H_t$ is a subgroup of the additive group of $k$ with $| H_t | = p^t$ and $b_j \in k^*$.
\end{lem}

\begin{table}
\begin{center}
\begin{tabular}{cccc}
$Case$ & $\G$ & $z$ & $Ramification$  \\
\hline \\
1 & $C_m$, $(m,p)=1$& $x^m$ & $(m,m)$\\ \\
2 & $D_{2m}$, $(m,p)=1$& $x^m+\frac{1}{x^m}$ & $(2,2,m)$\\ \\
3 & $A_4, \, p\neq 2, 3$ & $\frac{x^{12}-33x^8-33x^4+1}{x^2(x^4-1)^2}$ & $(2,3,3)$\\ \\
4 & $S_4, \, p\neq 2, 3$ & $\frac{(x^8+14x^4+1)^3}{108(x(x^4-1))^4}$ & $(2,3,4)$\\ \\
5 & $A_5, \,p\neq 2, 3, 5$ & $\frac{(-x^{20}+228x^{15}-494x^{10}-228x^5-1)^3}{(x(x^{10}+11x^5-1))^5}$ & $(2,3,5)$\\ \\
  & $A_5, \,p=3$ & $\frac{(x^{10}-1)^6}{(x(x^{10}+2ix^5+1))^5}$ & $(6,5)$\\ \\
6 & $U$ & $  \displaystyle{\prod_{a \in H_t}} (x+a)$ & $(p^t)$\\ \\
7 & $K_m$ & $(x   \displaystyle{\prod_{j=0}^{\frac{p^t-1}{m}-1} } (x^m-b_j))^m$ & $(mp^t,m)$ \\ \\
8 & $\psl(2,q), \,p\neq 2$ & $\frac{((x^q-x)^{q-1}+1)^{\frac{q+1}{2}}}{(x^q-x)^{\frac{q(q-1)}{2}}}$ & $(\alpha,\beta)$ \\ \\
9 & $\pgl(2,q)$ & $\frac{((x^q-x)^{q-1}+1)^{q+1}}{(x^q-x)^{q(q-1)}}$ & $(2\alpha,2\beta)$ \\ \\
\end{tabular}
\caption{Rational functions corresponding to each $\G$} \label{t1}
\end{center}
\end{table}

\def\bC{\textbf{C}}
\def\d{\delta}

Next we want to determine groups which occur as automorphism group $G$ of genus $g\geq 2$ cyclic curves, their signatures and the dimension of the corresponding locus. 
We know that $\G:=G/G_0$, where $G_0:=\Gal(k(x,y)/k(x))$ and  $\G$ is isomorphic to $C_m$, $D_m$, $A_4$, $S_4$, $A_5$, $U$, $K_m$, $\psl(2,q)$, $\pgl(2,q)$. By considering the lifting of ramified points in each $\G$, we divide each $\G$ into subcases. We determine the signature of each subcase by looking at the behavior of lifting and ramification of $\G$. Using the signature and \cref{e2} we calculate $\delta$ for each case. We list all possible automorphism groups $G$ as separate theorems for each $\G$.

We assume that   $5 < p \leq 2g+1$.   The case $p>2g+1$ is same as $ p =0$; see \cite{serdica}

\begin{thm} \label{th1}
Let $g \geq$ 2 be a fixed integer, $\X$ a genus $g$ cyclic curve, $G=\Aut(\X)$ and $C_n\triangleleft G$ such that $g(\X^{C_n})=0$. The signature of cover $\Phi : \X \to \X^G$ and dimension $\delta$ is given in \cite[Table~2]{sa-1}, where  $m=|\psl(2,q)|$ for cases 38-41 and $m=|\pgl(2,q)|$ for cases 42-45.
\end{thm}

There are 45 signatures from the above theorem (not all occur in every genus $g\geq 2$).   The following theorem gives us all possible automorphism groups of genus $g\geq 2$ cyclic curves defined over the finite field of characteristic $p$.

\begin{thm} \label{th14}
Let $\X_g$ be a genus $g\geq2$ irreducible cyclic curve defined over an algebraically closed field $k$, $\ch k=p\neq2$. Assume the cyclic group $C_n \lhd G=\Aut(\X_g)$ and let  $\G= \Aut(\X_g)/C_n $ be its reduced automorphism group.

\begin{enumerate}

\item If $\G \cong C_m$ then $G \cong C_{mn}$ or
\begin{center}
$\left\langle \r, \s \right|\r^n=1,\s^m=1,\s\r\s^{-1}=\r^l \rangle$
\end{center}
where (l,n)=1 and $l^m\equiv 1$ (mod n).

\item If $\G \cong D_{2m}$ then $G \cong D_{2m} \times C_n$ or
\begin{align*}
\begin{split}
G_{5}=& \left\langle \r, \s, \t \right|\r^n=1,\s^2=\r,\t^2=1,(\s\t)^m=1,\s\r\s^{-1}=\r,\t\r\t^{-1}=\r^{n-1} \rangle\\
G_{6}=& D_{2mn} \\
G_{7}=& \left\langle \r, \s, \t \right|\r^n=1,\s^2=\r,\t^2=\r^{n-1},(\s\t)^m=1,\s\r\s^{-1}=\r,\t\r\t^{-1}=\r \rangle\\
G_{8}=& \left\langle \r, \s, \t \right|\r^n=1,\s^2=\r,\t^2=1,(\s\t)^m=\r^{\frac{n}{2}},\s\r\s^{-1}=\r,\t\r\t^{-1}=\r^{n-1} \rangle \\
G_{9}=& \left\langle \r, \s, \t \right|\r^n=1,\s^2=\r,\t^2=\r^{n-1},(\s\t)^m=\r^{\frac{n}{2}},\s\r\s^{-1}=\r,\t\r\t^{-1}=\r \rangle\end{split}
\end{align*}

\item If $\G \cong A_4$ and $p\neq3$ then $G \cong A_4 \times C_n$ or
\begin{align*}
\begin{split}
G'_{10}=& \left\langle \r, \s, \t \right|\r^n=1,\s^2=1,\t^3=1,(\s\t)^3=1,\s\r\s^{-1}=\r,\t\r\t^{-1}=\r^l \rangle\\
G'_{12}=& \left\langle \r, \s, \t \right|\r^n=1,\s^2=1,\t^3=\r^{\frac{n}{3}},(\s\t)^3=\r^{\frac{n}{3}},\s\r\s^{-1}=\r,\t\r\t^{-1}=\r^l \rangle\\
\end{split}
\end{align*}
where $(l,n)=1$ and $l^3\equiv 1$ (mod n) or
\begin{center}
$\left\langle \r, \s, \t \right|\r^n=1,\s^2=\r^{\frac{n}{2}},\t^3=\r^{\frac{n}{2}},(\s\t)^5=\r^{\frac{n}{2}},\s\r\s^{-1}=\r,\t\r\t^{-1}=\r \rangle $
\end{center}
or
\begin{align*}
\begin{split}
G_{10}=& \left\langle \r, \s, \t \right|\r^n=1,\s^2=1,\t^3=1,(\s\t)^3=1,\s\r\s^{-1}=\r,\t\r\t^{-1}=\r^k \rangle\\
G_{13}=& \left\langle \r, \s, \t \right|\r^n=1,\s^2=\r^{\frac{n}{2}},\t^3=1,(\s\t)^3=1,\s\r\s^{-1}=\r,\t\r\t^{-1}=\r^k \rangle\\
\end{split}
\end{align*}
where $(k,n)=1$ and $k^3\equiv 1$ (mod n).

\item If $\G \cong S_4$ and $p\neq3$ then $G \cong S_4 \times C_n$ or
\begin{align*}
\begin{split}
G_{16}=& \left\langle \r, \s, \t \right|\r^n=1,\s^2=1,\t^3=1,(\s\t)^4=1,\s\r\s^{-1}=\r^l,\t\r\t^{-1}=\r \rangle\\
G_{18}=& \left\langle \r, \s, \t \right|\r^n=1,\s^2=1,\t^3=1,(\s\t)^4=\r^{\frac{n}{2}},\s\r\s^{-1}=\r^l,\t\r\t^{-1}=\r \rangle\\
G_{20}=& \left\langle \r, \s, \t \right|\r^n=1,\s^2=\r^{\frac{n}{2}},\t^3=1,(\s\t)^4=1,\s\r\s^{-1}=\r^l,\t\r\t^{-1}=\r \rangle \\
G_{22}=& \left\langle \r, \s, \t \right|\r^n=1,\s^2=\r^{\frac{n}{2}},\t^3=1,(\s\t)^4=\r^{\frac{n}{2}},\s\r\s^{-1}=\r^l,\t\r\t^{-1}=\r \rangle\\
\end{split}
\end{align*}
where $(l,n)=1$ and $l^2\equiv 1$ (mod n).

\item If $\G \cong A_5$ and $p\neq5$ then $G \cong A_{5}\times C_{n}$ or
\begin{center}
$\left\langle \r, \s, \t \right|\r^n=1,\s^2=\r^{\frac{n}{2}},\t^3=\r^{\frac{n}{2}},(\s\t)^5=\r^{\frac{n}{2}},\s\r\s^{-1}=\r,\t\r\t^{-1}=\r \rangle $
\end{center}

\item If $\G \cong U$ then $G \cong U \times C_n$ or
\begin{multline*}
<\r,\s_1,\s_2,...,\s_t|\r^n=\s_1^p=\s_2^p=...=\s_t^p=1,\\ \s_i\s_j=\s_j\s_i, \s_i\r\s_i^{-1}=\r^{l}, 1\leq i,j\leq t>
\end{multline*}
where $(l,n)=1$ and $l^p \equiv 1$ (mod n).

\item If $\G \cong K_m$ then $G \cong$

\begin{multline*}
<\r,\s_1,...,\s_t,v|\r^n=\s_1^p=...=\s_t^p=v^m=1, \s_i\s_j=\s_j\s_i,\\
v\r v^{-1}=\r, \s_i\r\s_i^{-1}=\r^{l}, \s_iv\s_i^{-1}=v^{k}, 1\leq i,j\leq t >
\end{multline*}
where $(l,n)=1$ and $l^p \equiv 1$ (mod n), $(k,m)=1$ and $k^p \equiv 1$ (mod m) or
\begin{align*}
\begin{split}
\left\langle \r,\s_1,...,\s_t|\r^{nm}=\s_1^p=...=\s_t^p=1, \s_i\s_j=\s_j\s_i, \s_i\r\s_i^{-1}=\r^{l}, 1\leq i,j\leq t\right\rangle
\end{split}
\end{align*}
where $(l,nm)=1$ and $l^p \equiv 1$ (mod nm).

\item If $\G \cong \psl(2,q)$ then $G\cong \psl(2,q)\times C_n$ or $SL_2(3)$.

\item If $\G \cong \pgl(2,q)$ then $G \cong \pgl(2,q) \times C_n$.

\end{enumerate}
\end{thm}

In \cite{s-sh} the corresponding equations are   given for each case.  In \cite{h-sh} for each group, it is discussed if the corresponding curve is defined over its field of moduli.


Applying \cref{th14}  we obtain the automorphism groups of all genus 3 superelliptic curves  defined over algebraically closed field of characteristic $p \neq 2$. Below we list the GAP group ID's of each of those groups.

\begin{lem}\label{thm_g_3}
Let $\X_g$ be a genus 3 superelliptic curve  defined over a field of characteristic $p \neq 2$. Then the automorphism groups of $\X_g$ is one of the following.

\begin{description}

    \item[i)] $p=3$ : $(2,1)$, $(4,2)$, $(3,1)$, $(4,1)$, $(8,2)$, $(8,3)$, $(7,1)$, $(14,2)$, $(6,2)$, $(8,1)$, $(8,5)$, $(16,11)$, $(16,10)$, $(32,9)$, $(30,2)$, $(16,7)$, $(16,8)$, $(6,2)$.

    \item[ii)] $p=5$ : $(2,1)$, $(4,2)$, $(3,1)$, $(4,1)$, $(8,2)$, $(8,3)$, $(7,1)$, $(21,1)$, $(14,2)$, $(6,2)$, $(12,2)$,
    $(9,1)$, $(8,1)$, $(8,5)$, $(16,11)$, $(16,10)$, $(32,9)$, $(42,3)$, $(12,4)$, $(16,7)$, $(24,5)$, $(18,3)$,
    $(16,8)$, $(48,33)$, $(48,48)$.

    \item[iii)] $p=7$ : $(2,1)$, $(4,2)$, $(3,1)$, $(4,1)$, $(8,2)$, $(8,3)$, $(7,1)$, $(21,1)$, $(6,2)$, $(12,2)$,
    $(9,1)$, $(8,1)$, $(8,5)$, $(16,11)$, $(16,10)$, $(32,9)$, $(30,2)$, $(42,3)$, $(12,4)$, $(16,7)$, $(24,5)$, $(18,3)$,
    $(16,8)$, $(48,33)$, $(48,48)$.

    \item[iv)] $p=0$ or $p \geq 11$ : $(2,1)$, $(4,2)$, $(3,1)$, $(4,1)$, $(8,2)$, $(8,3)$, $(7,1)$, $(21,1)$, $(14,2)$, $(6,2)$, $(12,2)$,
    $(9,1)$, $(8,1)$, $(8,5)$, $(16,11)$, $(16,10)$, $(32,9)$, $(30,2)$, $(42,3)$, $(12,4)$, $(16,7)$, $(24,5)$, $(18,3)$,
    $(16,8)$, $(48,33)$, $(48,48)$.
\end{description}

\end{lem}

Recall that the list for  $p=0$ is the same as for $p > 7$.   



Again applying \cref{th14}, we obtain the possible automorphism groups of genus 4 cyclic curves defined over an algebraically closed field of characteristic 0,3,5,7 and bigger than 7. We list the GAP group ID of these groups in following theorem.

\begin{lem}
Let $\X_g$ be a genus 4 cyclic curve defined over a field of characteristic $p$. Then the automorphism group of $\X_g$ is one of the following.

\begin{description}

    \item[i)] $p=3$ : $(2,1)$, $(4,2)$, $(3,1)$, $(6,2)$, $(5,1)$, $(10,2)$, $(20,1)$, $(9,1)$, $(18,2)$, $(15,1)$, $(4,1)$, $(20,4)$, $(8,3)$, $(40,8)$, $(12,5)$, $(16,7)$, $(20,5)$, $(32,19)$, $(24,10)$, $(8,4)$, $(9,2)$, $(18,5)$.

    \item[ii)] $p=5$ : $(2,1)$, $(4,2)$, $(3,1)$, $(6,2)$, $(9,2)$, $(5,1)$, $(10,2)$, $(20,1)$, $(9,1)$, $(27,4)$, $(18,2)$, $(4,1)$, $(18,3)$, $(8,3)$, $(12,5)$, $(36,12)$, $(54,4)$, $(16,7)$, $(20,5)$, $(32,19)$, $(24,10)$, $(8,4)$, $(60,9)$, $(36,11)$, $(24,3)$, $(72,42)$, $(10,2)$, $(18,5)$.

    \item[iii)] $p=7$ : $(2,1)$, $(4,2)$, $(3,1)$, $(6,2)$, $(9,2)$, $(5,1)$, $(10,2)$, $(20,1)$, $(9,1)$, $(27,4)$, $(18,2)$, $(15,1)$, $(4,1)$, $(20,4)$, $(18,3)$, $(8,3)$, $(40,8)$, $(12,5)$, $(36,12)$, $(54,4)$, $(16,7)$, $(20,5)$, $(32,19)$, $(24,10)$, $(8,4)$, $(60,9)$, $(36,11)$, $(24,3)$, $(72,42)$.

    \item[iv)] $p=0$ or $p \geq 11$ : $(2,1)$, $(4,2)$, $(3,1)$, $(6,2)$, $(9,2)$, $(5,1)$, $(10,2)$, $(20,1)$, $(9,1)$, $(27,4)$, $(18,2)$, $(15,1)$, $(4,1)$, $(20,4)$, $(18,3)$, $(8,3)$, $(40,8)$, $(12,5)$, $(36,12)$, $(54,4)$, $(16,7)$, $(20,5)$, $(32,19)$, $(24,10)$, $(8,4)$, $(60,9)$, $(36,11)$, $(24,3)$, $(72,42)$.
\end{description}
\end{lem}

The above two lemmas are a simple search going through all the cases of the theorem, but they illustrate the idea that for cyclic (superelliptic curves) all isomorphism classes of curves can be written out (including a parametric equation in each case).

There is one case missing from all the results of this section, namely $p=2$.  Next we will try to explain that case.
%
%


A hyperelliptic curve over an algebraically closed field $k$ of characteristic 2 admits an Artin-Schreier model
\[ \X :y^2+y=R(x)\]
where $R(x)$ is a rational function with no poles of even order. An isomorphism to another curve
$ \X^\prime : y^2 + y = Q(x)$
induces an automorphism of the projective $x$-line.    One can then determine possible normal forms for $R(x)$, work out which fractional linear transformations preserve each form and see how these lift to $\X$ and interact with the hyperelliptic involution.

Now the details. Let $\X^\prime $ be a hyperelliptic curve of genus $g$ over an algebraically closed field $k$ of characteristic 2. We use an Artin-Schreier model ${y^\prime}^2 + y^\prime  =  R'(x)$. As a consequence of Hasse's theory we
can find a rational function $Q(x) \in k(x)$ such that the rational function $R'(x) + Q(x) + {Q(x)}^2$ has no poles of even order. Let $R(x) = R'(x) + Q(x) + {Q(x)}^2$ and $y=y^\prime+Q(x)$ to get a curve $\X$ in  normalized form $y^2 + y = R(x)$. Then $y$ is unique up to transformations of the form $y \longmapsto y + B(x)$, where $B(x)$ is a rational function of $x$.

Now, take two hyperelliptic curves, $\X : y^2 + y = R(x)$ and $\X^\prime : y^2 + y = Q(x)$. Then, given  an isomorphism  $\phi : \X \longmapsto \X^\prime$ and the finite morphisms $f_1 : \X \longmapsto \P^{1}$, and $f_2 : \X^\prime \longmapsto \P^{1}$ of degree 2, there exists a unique automorphism $\sigma$ of $\P^{1}$ such that $f_{2}\circ \phi = \sigma \circ f_{1}$. Any isomorphism between these curves has the form
\[ (x, y) \longmapsto (\sigma(x),y+S(x)) = \left( \frac{ax + b}{cx + d}, \,  y + S(x) \right) \]
for some $S(x) \in k(x)$. Hence, these curves are isomorphic if and only if
\[ Q(x) = R \left( \frac{ax + b}{cx + d} \right) + S(x) + {S(x)}^2. \]

Let $\rm{div}(R)_{\infty} = \Sigma n_a(a)$ be the polar divisor of $R(x)$ on the projective line $\P^{1}$. $\X$ is ramified at each $a$ and if $P_a$ is the unique point of $\X$ over $a$ then the curve $y^2 + y = R(x)$ has the different
\[ \diff ( \X /\P^{1}) = \Sigma (n_a + 1) P_a, \]
where the $n_a$ are odd (\cite{Sti}, Prop III.$7.8$)
\[ 2g - 2 = -2 [F: k(x)] + \degrm ( \diff (\X/\P^{1})) \implies \degrm( \diff (\X /\P^{1})) = 2g + 2\]
The ramification types determine the isomorphism classes of the hyperelliptic curves. The solutions of the equation $\Sigma (n_a + 1) = 2g + 2$ in the unknown odd positive integers give us the
following ramification types:  $ (1, 1, 1, 1)$, $(3, 1, 1)$, $(3, 3)$, $(5, 1)$, $(7)$   for   genus  $g=3$.
Therefore we get the following normal forms.
\begin{equation}
\begin{split}
& y^2 + y = \left\{ \begin{array}{l} {\alpha}_1 x + {\alpha}_2 x^{-1} + {\alpha}_3 (x - 1)^{-1} +
{\alpha}_4
(x - \lambda)^{-1}\\
x^3 + \alpha x + \beta x^{-1} + \gamma (x - 1)^{-1}\\
x^3 + \alpha x + \beta x^{-3} + \gamma x^{-1}\\
x^5 + \alpha x^3 + \beta x^{-1}\\
x^7 + \alpha x^5 + \beta x^3
\end{array} \right.\\
\end{split}
\end{equation}

For genus $g=4$ we have
$(1, 1, 1, 1, 1)$, $(3, 1, 1, 1)$, $(3, 3, 1)$, $(5, 1, 1)$, $(5, 3)$, $(7, 1)$, $(9) $     and
therefore we get the following normal forms for genus $3$ and $4$ respectively.

\begin{equation}
\begin{split}
&  y^2 + y = \left\{ \begin{array}{l} {\alpha}_1 x + {\alpha}_2 x^{-1} + {\alpha}_3 (x - 1)^{-1} +
{\alpha}_4
(x - \lambda)^{-1} + {\alpha}_5 (x - \mu)^{-1}\\
x^3 + \alpha x + {\beta}_1 x^{-1} + {\beta}_2 (x - 1)^{-1} + {\beta}_3
(x - \lambda)^{-1}\\
x^3 + \alpha x + \beta x^{-3} + \gamma x^{-1} + \sigma (x - 1)^{-1}\\
x^5 + \alpha x^3 + \beta x^{-1} + \gamma (x - 1)^{-1}\\
x^5 + \alpha x^3 + \beta x^{-3} + \gamma x^{-1}\\
x^7 + \alpha x^5 + \beta x^3 + \gamma x^{-1}\\
x^9 + {\alpha}_1 x^7 + {\alpha}_2 x^5 + {\alpha}_3 x^3
\end{array} \right.
\end{split}
\end{equation}

Using a case by case analysis in \cite{demir} it is proved that

\begin{thm}
Let $C$ be a genus $g\geq 2$    hyperelliptic curve defined over a field $k$ such that $\ch k =2$.   Then

i) if $g=3$ then $\Aut(\X)$ is isomorphic to $C_2$, $C_4$, $V_4$, $C_2^3$, $C_6$, $C_{14}$, $D_{12}$

ii) if $g=4$ then $\Aut(\X)$ is isomorphic to $C_2$, $V_4$, $C_4$, $C_2^3$, $C_6$,  $C_{18}$, $D_{20}$.
\end{thm}

The corresponding equations are given in each case.  The higher genus cases can be determined in similar way.

\begin{rem} It seems as the above methods can be extended to determine complete lists of automorphism groups of any superelliptic curve of genus $g\geq 2$ and any $2  \leq p \leq 2g+1$.
We are not aware if such lists  are determined for $p=2, 3, 5$ for large $g$.
\end{rem}

We still have to consider the case $p=0$ and equivalently $p> 2g+1$.  However, here we can use the theory of compact Riemann surfaces and Fuchsian groups to give a complete answer to the question of determining the list of automorphism groups for any $g\geq 2$.

\begin{rem}Isomorphism classes of superelliptic curves are determined by invariants of binary forms.  For superelliptic curves with extra automorphisms other invariants are introduced in
\cite{g-sh}, \cite{antoniadis}. Finding equations of such superelliptic curves over a field of moduli is interesting on its own right; see \cite{h-sh}.
\end{rem}


\section{Automorphism groups of compact Riemann surfaces}\label{sec-CRS}

To contrast with the case $p>0$, and to provide some historical footing, we now turn our attention to the classical case: compact Riemann surfaces. Due to a contemporary adaptation of Riemann's Existence Theorem, see \cref{th-riem}, and its amenability to computational group theory, in modern work one of the most utilized tools for classification of automorphism groups of compact Riemann surfaces is also one of the more classical ones: uniformization and Fuchsian groups. In this section, we provide a short exposition outlining this traditional approach, describe some of the current results and how, with a little further direction, this method can be used to help determine full automorphism groups.

\subsection{Fuchsian Groups and Signatures}\label{subsec-FGS}

The Uniformization Theorem states that any compact Riemann surface $\X$ of genus $g \geq 2$ is conformally equivalent to a quotient of its universal cover, the upper half plane ${\mathbb H}$, by a torsion free discrete subgroup $\Lambda$ of $\Aut({\mathbb H}) =\psl(2,{\mathbb R})$. The group $\Lambda$ is isomorphic to $\Pi_{g}$, the fundamental group of $\X$, and is called a {\em surface group} for $\X$. For a given surface $\X$, surface groups are unique up to conjugacy in $\psl(2,{\mathbb R})$ meaning two compact Riemann surfaces $\X$ and $\X'$ are conformally equivalent if and only if surface groups for $\X$ and $\X'$ are conjugate in $\psl(2,{\mathbb R})$.

Now if $\phi \colon \X\rightarrow \X$ is an automorphism of $\X$, then it can be lifted to an automorphism $\phi_{\Lambda}$ of ${\mathbb H}$ which normalizes $\Lambda$. In particular, if $G$ is a group of automorphisms of $\X$, then $G$ can be lifted to a discrete subgroup $\Gamma$ of $\psl(2,{\mathbb R})$, called a {\em Fuchsian group}, containing $\Lambda$ with index $|G|$ and which normalizes $\Lambda$. We call $\Gamma$ the Fuchsian group corresponding to $G$, and if $\Lambda$ has been fixed, we call $G$ the automorphism group corresponding to $\Gamma$.

Conversely, if $\Gamma$ is a Fuchsian group and $\Lambda$ is a normal subgroup of $\Gamma$ which is isomorphic to $\Pi_{g}$, then there is a natural action of the quotient group $\Gamma /\Lambda$ on the quotient surface ${\mathbb H} /\Lambda$ which is a surface of genus $g$.

These observations illustrate the primary basic tool for determining group actions on compact Riemann surfaces of genus $g \geq 2$: for a given genus, determine all Fuchsian groups (up to isomorphism) for which there exists a normal subgroup isomorphic to $\Pi_{g}$. To explain this process in more detail, we need some additional preliminary results.

For a cocompact Fuchsian group $\Gamma$, the quotient surface ${\mathbb H}/\Gamma$ is a compact Riemann surface and the quotient map $\pi_{\Gamma} \colon {\mathbb H} \rightarrow {\mathbb H}/\Gamma$ is a conformal, possibly branched, map. We define the {\em signature} of $\Gamma$ to be the tuple $(g_{\Gamma} ;m_{1}, m_{2},\dots ,m_{r})$ where $g_{\Gamma}$ is the genus of ${\mathbb H} /\Gamma$ and the quotient map  $\pi_{\Gamma}$ branches over $r$ points with ramification indices $m_{i}$ for $1\leq i\leq r$. We call $g_{\Gamma}$ the {\bf{orbit genus}} of $\Gamma$ and the numbers $m_{1},\dots ,m_{r}$ the {\bf{periods}} of $\Gamma$. The signature of $\Gamma$ provides information regarding a presentation for $\Gamma$:

\begin{thm}
\label{Sig}
If $\Gamma$ is a Fuchsian group with signature $(g_{\Gamma};m_{1},\dots ,m_{r})$ then there exist group elements $\alpha_{1}, \beta_{1} ,\dots ,\alpha_{g_{\Gamma}} , \beta_{g_{\Gamma}},\gamma_{1},\dots \gamma_{r}\in \psl(2,{\mathbb R})$, such that;

\begin{enumerate}
\item
$\Gamma=\langle \alpha_{1},\beta_{1},\dots ,\alpha_{g_{\Gamma}},\beta_{g_{\Gamma}},\gamma_{1},\dots \gamma_{r}\rangle $.

\item
Defining relations for $\Gamma$ are

\begin{center}
$\gamma_{1}^{m_{1}}, \gamma_{2}^{m_{2}},\dots ,\gamma_{r}^{m_{r}},\prod\limits_{i=1}^{g_{\Gamma}} [\alpha_{i},\beta_{i}] \prod\limits_{j=1}^{r} \gamma_{j}$.
\end{center}

\item
Each elliptic element (the elements of finite order) lies in a unique conjugate of $\langle \gamma_{i}\rangle $ for suitable $i$. Furthermore, the cyclic groups $\langle \gamma_{i}\rangle $ are self-normalizing in $\Gamma$.

\item
Each elliptic element of $\Gamma$ has a unique fixed point in ${\mathbb H}$. All other elements (the hyperbolic elements) act fixed point freely on ${\mathbb H}$.

\end{enumerate}

\end{thm}

\noindent
We call a set of elements of $\Gamma$ satisfying \cref{Sig} {\bf{canonical generators}} for $\Gamma$. Notice that if $\Gamma$ is a surface group for a surface of genus $g$, since it is torsion free, it must have signature $(g;-)$.

Now, if $G$ acts conformally on $\X$, $\Lambda$ is surface group for $\X$, and $\Gamma$ is the Fuchsian group corresponding to $G$ with signature $(g_{\Gamma};m_{1},\dots ,m_{r})$, then there exists an epimorphism, called a {\bf surface kernel epimorphism}, $\rho \colon \Gamma \rightarrow G$ with kernel $\Lambda$. This epimorphism can be neatly summarized in the context of finite groups by a {\bf generating vector}. Specifically, if $\alpha_{1}, \beta_{1} ,\dots ,\alpha_{g_{\Gamma}} , \beta_{g_{\Gamma}},\gamma_{1},\dots \gamma_{r}\in \psl(2,{\mathbb R})$ are canonical generators for $\Gamma$, then we get a $\left(  2g_{\Gamma}+r\right)  $-tuple of elements from $G$,
$(a_{1},\ldots a_{g_{\Gamma}},b_{1},\ldots b_{g_{\Gamma}},c_{1},\ldots c_{r})$ where $\rho (\alpha_i ) =a_i$, $\rho (\beta_i )=b_i$ and $\rho (\gamma_i)=c_i$ called a
$(g_{\Gamma} ;m_{1},\ldots,m_{r})$-generating vector for $G$. Moreover, since $\rho$ is an epimorphism with torsion free kernel, we have:
\begin{itemize}
\item
$\displaystyle \prod\limits_{i=1}^{g_{\Gamma}} [a_{i},b_{i}] \prod\limits_{j=1}^{r} c_{j}=1$ (the identity)
\item
${\rm O}(c_i)=m_{i}$ where ${\rm O}$ denotes element order
\end{itemize}
Moreover, provided a set of canonical generators of $\Gamma$ have been fixed, there is a one to one correspondence between the set of $(g_{\Gamma} ;m_{1},\ldots,m_{r})$ - generating vectors of $G$ and $\rm{Epi} (\Gamma,G),$ epimorphisms $\Gamma\rightarrow G$ preserving the orders of the $\gamma_{j}$.

Now suppose that $\Gamma$ is a Fuchsian group with signature $(g_{\Gamma};m_{1},\dots ,m_{r})$ and suppose that $\Lambda$ is a normal surface subgroup of $\Gamma$ for a surface of genus $g$. Letting $G=\Gamma/\Lambda$ and identifying the orbit spaces ${\mathbb H} /\Gamma$ and $\X/G$ we get the tower of covers and quotient maps given in \cref{Quot}.

\begin{figure}[h]
$
\xymatrix @R=.5in @C=1in {
{\mathbb H} \ar@/^2pc/[rr]^{\pi_{\Gamma}} \ar[r]_{\pi_{\Lambda}} & {\mathbb H} /\Lambda =\X \ar[r]_{\pi_{G}} & {\mathbb H} /\Gamma  =\X/G \\
}
$
\caption{\label{Quot}Holomorphic Quotient Maps and Surface Identifications}
\end{figure}

Since the universal covering map $\pi_{\Lambda}$ is unramified, it follows that the quotient map $\pi_{G} \colon \X\rightarrow \X/G$ is branched over the same points as $\pi_{\Gamma} \colon {\mathbb H} \rightarrow {\mathbb H} /\Gamma$ with the same ramification indices. In particular, $\X/G$ has genus $g_{\Gamma}$ and $\pi_{G}$ is a degree $|G|$ map branched over $r$ points with ramification indices $m_1,\dots m_r$. Consequently,  we often say that $G$ acts or has the same signature as $\Gamma$.

 Now since the map $\pi_{G}$ is between compact Riemann surfaces, the Riemann-Hurwitz formula holds, giving: $$g  -1=|G|(g_{\Gamma}-1) +\frac{|G|}{2} \sum_{j=1}^{r} \bigg( 1-\frac{1}{m_{j}} \bigg).$$ Combining our observations, we get the following modern adaption of Riemann's existence theorem which provides necessary and sufficient conditions for the existence of the action of a group $G$ on a compact Riemann surface of genus $g$ with signature $(g_{\Gamma};m_1,\dots ,m_r)$:

\begin{thm}
\label{th-riem}
A finite group $G$ acts on a compact Riemann surface $S$ of genus $g \geq 2$ with signature $(g_{\Gamma};m_{1},\dots ,m_{r})$ if and only if:
\begin{enumerate}
\item
the Riemann--Hurwitz formula is satisfied:  $$g -1=|G|(g_{\Gamma}-1) +\frac{|G|}{2} \sum_{j=1}^{r} \bigg( 1-\frac{1}{m_{j}} \bigg).$$
\item
there exists an $(g_{\Gamma};m_{1},\dots ,m_{r})$-generating vector for $G$.
\end{enumerate}
\end{thm}


\begin{rem}\label{rk-fundgroup}
Let $\Y$ denote the orbit surface $\X/G=\mathbb{H}/\Gamma$ and $\pi_{G}:\X\rightarrow\Y$ the
quotient morphism. The set of branch points over which $\pi_{G}$ is ramified is denoted $\mathcal{B}=\mathcal{B}_{G}.$ The covering $\pi_{G}:\X^{\circ}\rightarrow\Y^{\circ},$ where $\Y^{\circ
}=\Y-\mathcal{B}_{G}$ and $\X^{\circ}=\X-\pi
_{G}^{-1}(\mathcal{B}_{G})$ is an unramified Galois covering of affine curves.
The fundamental group $\pi_{1}(\Y^{\circ})$ has a presentation
\[
\pi_{1}(\Y^{\circ})=\langle\alpha_{1},\beta_{1},\dots,\alpha_{g_{\Gamma}},\beta_{g_{\Gamma}},\gamma_{1},\dots\gamma_{r}:\prod \limits_{i=1}^{g_{\Gamma}}[\alpha_{i},\beta_{i}]\prod\limits_{j=1}^{r}\gamma_{j}=1\rangle.
\]
Any Galois cover $\X\rightarrow\Y$ , ramified exactly over
$\mathcal{B}_{G}$, with Galois group $G$, is defined by an epimorphism
$\eta:\pi_{1}(\Y^{\circ})\rightarrow G,$ and, hence, defined by a
generating vector, where we do not impose specific orders on the $c_{j}.$ The
generating vectors classify (not uniquely) the finite covers of $\Y$
, branched over $\mathcal{B}.$ There are only finitely many covers for a given
pair ($\Y,\mathcal{B)}$ and group order $\left\vert G\right\vert .$
This surface construction approach allows us adopt a similar, though unwieldy,
approach for $p>0$ using the \'{e}tale fundamental group $\pi_{1}^{et}(\Y^{\circ})$, discussed in \cref{rk-algfundgroup}
\end{rem}


\subsection{Translating the Problem into Finite Group Theory}

The importance of \cref{th-riem} is that it translates the problem of determining group actions of compact Riemann surfaces from a problem about infinite discrete groups into a problem about finite groups through the introduction of generating vectors. In particular, it makes the problem amenable to computational group theory and accordingly, classification results have significantly improved over the last few decades.

There are many different approaches to classifying automorphism groups, but perhaps the most common approach is to do so by genus; that is, fix a genus $g$ and then find all possible automorphism groups that can act on a surface of that genus and the signatures with which they act. The basic approach to this form of classification for a fixed genus $g$ is as follows:

\begin{enumerate}
\item
Find all possible signatures for each possible group order which satisfy the Riemann-Hurwitz formula.

\item
For each group order and signature, run over all groups of that possible order to either construct a generating vector, or show that no such generating vector exists.

\end{enumerate}

We illustrate with a couple of examples.

\begin{exa}
\label{ex-genus2}
Suppose that $g =2$ and consider the group order $|G|=2$. Then the Riemann-Hurwitz formula gives: $$2 -1=1=2(g_{\Gamma}-1) +\frac{2}{2} \sum_{j=1}^{r} \bigg( 1-\frac{1}{m_{j}} \bigg).$$ Since $|G|=2$, and the $m_i$'s are all element orders of $G$, it follows that $m_i=2$ for all $i$, so simplifying, we get $$3=2g_{\Gamma}+\frac{r}{2}.$$ Solving, we get $g_{\Gamma}=0$ and $r=6$ or $g_{\Gamma}=1$ and $r=2$, so the signatures $(0;2,2,2,2,2,2)$ and $(1;2,2)$. Now the only group $G$ of order $2$ is cyclic, so let $x$ be a generator of a group $G$ of order $2$. Then $(x,x,x,x,x,x)$ is a $(0;2,2,2,2,2,2)$ and $(x,x,x,x)$ is a $(1;2,2)$-generating vector for $G$, and hence $G$ acts on a surface of genus $2$ with signatures $(0;2,2,2,2,2,2)$ and $(1;2,2)$ but no other signatures.
\end{exa}

As illustrated in \cref{ex-genus2}, for a fixed group order and genus $g$, after simplification we obtain a linear equality in $k+1$ variables where $k$ is the number of distinct element orders in $G$. However, since all the variables in the equality are non-negative integers, since we have fixed $g$ and the group order $|G|$, there are just finitely many solutions and hence the problem can be solved fairly easily by computer (or even by hand). Of course, complete classification relies on running over all possible group orders in a given genus, but from the Hurwitz bound, we know there are only finitely many possibilities. In particular, the problem of classification can be solved completely computationally through an algorithm similar to the following:

\begin{alg}
\label{alg-class}
For $2\leq n\leq 84(g -1)$, we do the following:

\begin{enumerate}
\item
Solve for all signatures satisfying $$g -1=n(g_{\Gamma}-1) +\frac{n}{2} \sum_{j=1}^{r} \bigg( 1-\frac{1}{m_{j}} \bigg)$$ where the $m_j$ are divisors of $n$.

\item
For each signature $(g_{\Gamma};m_1,\dots m_r)$ found in (a), we do the following:

\begin{enumerate}
\item
For each group $G$ of order $n$ with elements of orders $m_1,\dots ,m_r$, construct all vectors of elements of $G$ of length $2g_{\Gamma}+r$ where the first $2g_{\Gamma}$ elements are any elements of $G$, and the $2g_{\Gamma}+i$th element has order $m_i$.

\item
For each vector $(a_1,b_1,\dots ,a_{g_{\Gamma}},b_{g_{\Gamma}},c_1,\dots ,c_r)$, test the relation $$\displaystyle \prod\limits_{i=1}^{g_{\Gamma}} [a_{i},b_{i}] \prod\limits_{j=1}^{r} c_{j}=1.$$ If there exists a vector for which this relation holds, we have constructed a generating vector and hence an action of $G$ exists with this signature. If no such vector satisfies this relation, $G$ does not act with this signature.

\end{enumerate}

\end{enumerate}
\end{alg}

Using essentially \cref{alg-class}, one can determine all possible signatures for a given genus $g\geq2$.
There is a huge amount of literature on this; see \cite{kyoto} and \cite{hkt} for a complete list of references. As previously mentioned, some of the most comprehensive results are from \cite{Br} where Breuer was able to determine all possible signatures of genus up to $48$ using GAP.

\subsection{Full automorphism groups.}

\subsubsection{Candidates for Non-Maximal Automorphism Groups and Singerman's List}
For each group and signature in most of the classification data available, all subgroups are also in the list. This raises the question how to pick out those groups that occur as the {\bf full automorphism group} of some surface of genus $g$. One way to answer this question is to translate the problem back into a problem about Fuchsian groups. Specifically, suppose that $G$ acts as a group of automorphisms on $\X$ with signature $(g_{\Gamma};m_1,\dots m_r)$ and that for the surface group $\Lambda$, $\Gamma$ is the Fuchsian group corresponding to $G$. Then $G$ is not the full automorphism group of $\X$ if and only if there exists a Fuchsian overgroup $\Gamma_1$ of $\Gamma$ which also contains $\Lambda$ as a normal subgroup. Therefore, we need to understand the subgroup and overgroup structure of Fuchsian groups based on signature, and in particular when, for a given Fuchsian group $\Gamma$, there is an overgroup $\Gamma_1$. This exact problem was solved by Singerman in \cite{Si}, and we summarize his results.

\begin{thm}
\label{th-sing}
If $\Gamma$ is a Fuchsian group whose signature does not appear in the second column of \cref{List}, then $\Gamma$ is isomorphic to a finitely-maximal Fuchsian group, that is, a group that is not contained with finite index in any other Fuchsian group. If the signature of $\Gamma$ does appear in the second column of \cref{List}, then $\Gamma$ is a subgroup of a Fuchsian group $\Gamma_1$ with signature  from the third column with finite index given in the last column.
\end{thm}

\begin{table}[h]
{\begin{tabular}{||c|l|l|l||}
\hline
Case & Signature $\Gamma$ & Signature $\Gamma_1$ & $[\Gamma_1 :\Gamma ]$ \\ [1.5ex] \hline \hline
N1 & $(2;-)$ & $(0;2,2,2,2,2,2)$ & $2$ \\[1ex] \hline
N2 & $(1;t,t)$ & $(0;2,2,2,2,t)$ & $2$ \\[1ex] \hline
N3 & $(1;t)$ & $(0;2,2,2,2t)$ & $2$ \\[1ex] \hline
N4 & $(0;t,t,t,t)$, $t\geq 3$ & $(0;2,2,2,t)$ & $4$ \\[1ex] \hline
N5 & $(0;t,t,u,u)$, $t+u\geq 5$ & $(0;2,2,t,u)$ & $2$ \\[1ex] \hline
N6 & $(0;t,t,t)$, $t\geq 4$ & $(0;3,3,t)$ & $3$ \\[1ex] \hline
N7 & $(0;t,t,t)$, $t\geq 4$ & $(0;2,3,2t)$ & $6$ \\[1ex] \hline
N8 & $(0;t,t,u)$, $t\geq 3$,$t+u\geq 7$ & $(0;2,t,2u)$ & 2 \\[1ex] \hline
T1 & $(0;7,7,7)$ & $(0;2,3,7)$ & $24$ \\[1ex] \hline
T2 & $(0;2,7,7)$ & $(0;2,3,7)$ & $9$ \\[1ex] \hline
T3 & $(0;3,3,7)$ & $(0;2,3,7)$ & $8$ \\[1ex] \hline
T4 & $(0;4,8,8)$ & $(0;2,3,8)$ & $12$ \\[1ex] \hline
T5 & $(0;3,8,8)$ & $(0;2,3,8)$ & $10$ \\[1ex] \hline
T6 & $(0;9,9,9)$ & $(0;2,3,9)$ & $12$ \\[1ex] \hline
T7 & $(0;4,4,5)$ & $(0;2,4,5)$ & 6 \\[1ex] \hline
T8 & $(0;n,4n,4n)$, $n\geq 2$ & $(0;2,3,4n)$ & $6$ \\[1ex] \hline
T9 & $(0;n,2n,2n)$, $n\geq 3$ & $(0;2,4,2n)$ & $4$ \\[1ex] \hline
T10 & $(0;3,n,3n)$, $n\geq 3$ & $(0;2,3,3n)$ & $4$ \\[1ex] \hline
T11 & $(0;2,n,2n)$, $n\geq 4$ & $(0;2,3,2n)$ & $3$ \\[1ex] \hline
\end{tabular}
\caption{Singerman's List.\label{List}}}
\end{table}

The importance of \cref{th-sing} is that it tells us that unless $G$ acts with one of the signatures in the second column of \cref{List}, then there always exists a Fuchsian group with that signature that is not contained with finite index in any other Fuchsian group $\Gamma_1$. In particular, the surface $\X=\Gamma /{\rm Ker} (\rho)$ where $\rho$ is the surface kernel epimorphism from $\rho$ to $G$ given by the generating vector of $G$ has $G$ as its full group of automorphisms.

Alternatively, there is a \emph{moduli dimension} argument in \cite{kyoto} which determines which signatures give full automorphism groups.  Moreover, from methods in \cite{kyoto} for any fixed $g\geq 2$ one can determine completely inclusion among the subloci of the moduli space $\mathcal M_g$ for all the groups.  See for example such diagrams for $g=3, 4$ in \cite{bsz}. It is noted in \cite{bsz} that the majority of cases come from  (cyclic)  superelliptic curves, and as previously discussed, such curves are well understood.

\subsubsection{Conditions for when a Group is not a Full Automorphism Group}

Now, just because a signature does appear in Singerman's list does not necessarily mean that a corresponding automorphism group is always contained in some larger automorphism group -- just that it might. By considering containments of the Fuchsian groups given in Singerman's list however, it is possible to determine necessary and sufficient conditions in terms of generating vectors or surface kernel epimorphisms for when a group does extend to some larger group. This was the primary goal of \cite{CondBuj} which we summarize in the following two theorems.

\begin{thm}
\label{Cond1}
Let $G$ be a finite group acting with a non-maximal and non-triangular Fuchsian signature on a compact Riemann surface $\X$ of genus $g$.

\begin{enumerate}
\item
Suppose $G$ acts with signature $(2; -)$ and has corresponding $(2; -)$-generating vector $(a_1,b_1,a_2,b_2)$. Then $G$ is contained in some larger group of automorphisms with corresponding signature $(0;2,2,2,2,2,2)$ if and only if the assignment $a_1 \rightarrow a_1^{-1}$, $b_1\rightarrow a_1b_1^{-1}a_1^{-1}$, $a_2\rightarrow (b_1^{-1}a_2b_2)a_2^{-1}(b_1^{-1}a_2b_2)^{-1}$ and $b_2\rightarrow (b_1^{-1}a_2)b_2^{-1} (b_1^{-1}a_2)^{-1}$ is an automorphism of $G$.

\item
Suppose $G$ acts with signature $(1; t,t)$ and has corresponding $(1; t, t)$-generating vector $(a_1,b_1,c_1,([a_1,b_1]x)^{-1})$. Then $G$ is contained in some larger group of automorphisms with corresponding signature $(0;2,2,2,2,t)$ if and only if the assignment $a_1\rightarrow a_1^{-1}$, $b_1\rightarrow b_1^{-1}$ and $c_1\rightarrow (a_1b_1)^{-1} c_1^{-1} (ba)$ is an automorphism of $G$.

\item
Suppose $G$ acts with signature $(1;t)$ and has corresponding $(1,t)$-generating vector $(a_1,b_1,[a_1,b_1]^{-1})$. Then $G$ is contained in some larger group of automorphisms with corresponding signature $(0;2,2,2,2t)$ if and only if the assignment $a_1\rightarrow a_1^{-1}$, $b_1\rightarrow b_1^{-1}$ is an of $G$.

\item
Suppose $G$ acts with signature $(0; t, t, u, u)$, where $t+u > 5$ and has corresponding $(0; t, t, u, u)$-generating vector $(c_1,c_2,c_3,c_4)$. Then $G$ is contained in some larger group of automorphisms with corresponding signature $(0;2,2,t,u)$ if and only if the assignment $c_1\rightarrow c_2$, $c_2\rightarrow c_1$, $c_3\rightarrow c_1^{-1}c_4c_1$ and $c_4\rightarrow c_2c_3c_2^{-1}$ is an automorphism of $G$.

\end{enumerate}

\end{thm}

\begin{thm}
\label{Cond2}
Let $G$ be a finite group acting on a compact Riemann surface $\X$ of genus $g$ with a triangular signature $(0; m_1,m_2,m_3)$ with corresponding $(0; m_1,m_2,m_3)$-generating vector $(c_1,c_2,c_3)$. Then $G$ is not the full automorphism group of $\X$ if and only if at least one of the following conditions is satisfied (up to permutation of the periods $m_1$, $m_2$, $m_3$).

\begin{enumerate}
\item
$G$ acts with signature $(0;t,t,t)$ where $t > 4$, and the assignment $c_1\rightarrow c_2$, $c_2\rightarrow c_3$ and $c_3\rightarrow c_1$ induces an automorphism of $G$.

\item
$G$ acts with signature $(0; t, t, u)$ where $t > 3$ and $t+u > 7$, and the assignment $c_1 \rightarrow c_2$, $c_2\rightarrow c_1$, $c_3\rightarrow  c_2c_3c_2^{-1}$ induces an automorphism of $G$.

\item
$G$ acts with signature $(0; 2, 7, 7)$, the conjugates of $c_2c_3^{-1}c_2c_1c_3^3$ generate a normal subgroup $K$ of index $56$ in $G$, and $G$ is extendable to a group $G'$ containing $G$ as a subgroup of index $9$ such that $G'$ is generated by $c_1$ and an element $\alpha$ which normalizes $K$ and satisfies $\alpha^3=1$, $(c_1\alpha )^7=1$, $c_2=(\alpha c_1 \alpha)^{-1} c_1\alpha (\alpha c_1 \alpha)$ and $c_3=(\alpha c_1 \alpha^{-1}) c_1\alpha (\alpha c_1 \alpha^{-1})$.

\item
$G$ acts with signature $(0; 3, 3, 7)$, the conjugates of $c_2c_1c_3^2$ generate a normal subgroup $K$ of index $21$ in $G$, and $G$ is extendable to a group $G'$ containing $G$ as a subgroup of index $8$ such that $G'$ is generated by $c_3$ and an element $\alpha$ which normalizes $K$ and satisfies $\alpha^3 = 1$, $(\alpha c_3)^2 = 1$, $c_1 = \alpha c_3^{-2} \alpha c_3^2 \alpha^{-1}$ and $c_2 = \alpha^{-1} c_3^2\alpha c_3^{-2} \alpha$.

\item
$G$ acts with signature $(0; 3, 8, 8)$, conjugates of $c_2^2c_1c_3^2$ and $c_3^{-1}c_2c_1^{-1} c_2^{-1} c_1c_2^{-1}$ generate a normal subgroup $K$ of index $72$ in $G$, and $G$ is extendable to a group $G'$ containing $G$ as a subgroup of index $10$ such that $G'$ is generated by $c_3$ and an element $\alpha$ which normalizes $K$ and satisfies $\alpha^3 = 1$, $(\alpha c_3)^2= 1$, $c_1 = \alpha c_3^{-2}\alpha c_3^2 \alpha^{-1}$ and $c_2=\alpha^{-1} c_3^2\alpha^{-1} c_3\alpha c_3^{-2} \alpha$.

\item
$G$ acts with signature $(0; 4, 4, 5)$, the conjugates of $c_1^{-1}c_2^{-1}c_3^2$ generate a normal subgroup $K$ of index $20$ in $G$, and $G$ is extendable to a group $G'$ containing $G$ as a subgroup of index $6$ such that $G'$ is generated by $c_3$ and an element $\alpha$ which normalises $K$ and satisfies $\alpha^4=1$, $(\alpha c_3)^2=1$, $c_1=\alpha^2c_3\alpha c_3^{-1} \alpha^2$ and $c_2=\alpha^{-1} c_3\alpha c_3^{-1} \alpha$

\item
$G$ acts with signature $(0; 3, n, 3n)$ where $n > 3$, the conjugates of $c_2$ generate a normal subgroup $K$ of index $3$ in $G$, and $G$ is extendable to a group $G'$ containing $G$ as a subgroup of index $4$ such that $G'$ is generated by $c_3$ and an element $\alpha$ which normalizes $K$ and satisfies $\alpha^2 = 1$, $(c_3\alpha )^3 = 1$, $c_1 = \alpha c_3 (c_3\alpha )^{-1} c_3^{-1} \alpha$ and $c_2 = \alpha c_3^3\alpha$

\item
$G$ acts with signature $(0; 2, n, 2n$) where $n > 4$, the conjugates of $c_2$ generate a normal subgroup $K$ of index $2$ in $G$, and $G$ is extendable to a group $G'$ containing $G$ as a subgroup of index $3$ such that $G'$ is generated by $c_3$ and an element $\alpha$ which normalizes $K$ and satisfies $\alpha^3=1$, $c_1=\alpha (\alpha c_3) \alpha^{-1}$ and $c_2=\alpha^{-1} c_3^2\alpha$.

\end{enumerate}

\end{thm}

We illustrate with an example.

\begin{exa}
The vector $(x,x,x^3)$ is a $(0;5,5,5)$-generating vector for the cyclic group $C_5=\langle x\rangle$, so applying the Riemann-Hurwitz formula, we get a $C_5$ action on a surface $\X$ of genus $g =2$ with signature $(0;5,5,5)$. Letting $(a,b,c)=(x,x,x^3)$, we see that  the assignment $c_1\rightarrow c_2$, $c_2\rightarrow c_1$ and $c_3\rightarrow c_2c_3c_1^{-1}$ induces an automorphism of $G$ (the trivial automorphism), so it follows by (2) of \cref{Cond2} that $C_5$ is not the full automorphism group of $\X$.
\end{exa}

\subsubsection{Finding Full Automorphism Groups Using Counting Methods\label{count}}

\cref{Cond1} and \cref{Cond2} provide explicit ways to determine whether or not a given group is the full automorphism group of a compact Riemann surface by looking at the corresponding generating vectors and possible the overgroup structures. In certain cases, we can avoid the actual explicit computation of a generating vector to determine full automorphism groups by using a different method that relies on counting epimorphisms instead. An advantage to this alternate method is that the process is somewhat iterative and uses the existing lists (such as Breuer's) -- something  \cref{Cond1} and \cref{Cond2} do not do. The main disadvantage to this alternate methods is that it does not work in all cases, and in some cases, we still need to return to generating vectors.

A Fuchsian group $\Gamma$ with a signature of the form $(0;m_1,m_2,m_3)$ is unique up to conjugation in $\psl(2,{\mathbb R})$. It follows that if we fix a Fuchsian group $\Gamma$ with signature $(0;m_1,m_2,m_3)$ and $\X$ is any surface on which a finite group $G$ acts with signature $(0;m_1,m_2,m_3)$, then there is a surface group $\Lambda$ for $\X$ which is normal in $\Gamma$ with $\Gamma /\Lambda =G$. In particular, the number of distinct surfaces with such a $G$-action will be equal to the number of
non-$\psl(2,{\mathbb R})$-conjugate torsion free normal subgroups of $\Gamma$ with quotient $G$. Of course, in general it is difficult to determine whether two subgroups of $\psl(2,{\mathbb R})$ are conjugate, but for triangle groups we have the following proved in \cite{GirWol}:
\begin{thm}
\label{thm-GirWol}
If the $\psl(2,{\mathbb R})$-conjugate surface groups $\Lambda$ and $\Lambda '$ are both normal subgroups of the triangle group $\Gamma$, then $\Lambda '=\alpha \Lambda \alpha^{-1}$ for some $\alpha \in N(\Gamma )$ or $\alpha \in N(N(\Lambda ))$ (where $N(.)$ denotes normalizer).
\end{thm}
In particular, two normal surface subgroups of a triangle group $\Gamma$ will be $\psl(2,{\mathbb R})$-conjugate if and only if they are conjugate within some triangle group $\Gamma_{1}$ containing $\Gamma$, and we know all possible such pairs of triangle groups from Singerman's list.  Thus, for a given triangle group $\Gamma$ with signature $(0;m_1,m_2,m_3)$ and finite group $G$, we can count the number of distinct surfaces (up to conformal equivalence) on which $G$ acts with signature $(0;m_1,m_2,m_3)$ by doing the following:
\begin{enumerate}
\item
Count the number of epimorphisms from $\Gamma$ onto $G$.
\item
Divide this number by $|\Aut(G)|$, the size of the automorphism group of $G$. This will give the number of distinct surface subgroups of $\Gamma$ with quotient $G$.
\item
For each triangle group $\Gamma_1$ containing $\Gamma$, sort the surface subgroups of $\Gamma$ with quotient $G$ into $\Gamma_1$-conjugacy classes.
\end{enumerate}
We discuss these steps in a little more detail, and explain how they can help with the problem of determining full automorphism groups.

First, at least currently, there is no straightforward way to calculate the number of epimorphisms from $\Gamma$ onto $G$ without constructing generating vectors. Fortunately, there is a way to use the characters of $G$ to count torsion free homomorphisms. Specifically the following is a consequence of the more general main result from  \cite{Jones1}:
\begin{thm}
\label{Jones}
Suppose that $\Gamma$ has signature $(0;m_{1},m_{2},m_{3})$ and denote by $L_{i}$ for $1\leq i\leq 3$ the union of all conjugacy classes of $G$ of elements of orders $m_{i}$ respectively. The number of homomorphisms from $\Delta$ to $G$ with torsion free kernel, denoted $|{\rm Hom} (\Delta ,G)|$, is given by:

\begin{equation}\label{Calc}
\frac{1}{|G|} \sum_{\chi}
\left\{ \chi (1)^{-1} [\sum_{x_{1}\in L_{1}} \chi (x_{1} )] [\sum_{x_{2}\in L_{2}} \chi (x_{2})] [\sum_{x_{3}\in L_{3}} \chi (x_{3})] \right\}
\end{equation}

\noindent
where $\chi$ runs over the irreducible characters of $G$.
\end{thm}
Of course, we need to be a little careful as some of the homomorphisms counted in \cref{Jones} may not be surjective onto $G$, but rather define an epimorphism onto a strict subgroup $H$ of $G$. However, if this is the case, the Riemann-Hurwitz formula ensures a corresponding action of $H$ on some surface of smaller genus. In particular, if we are working iteratively through a list, then we already know the number of epimorphisms of $\Gamma$ onto $H$. Thus to refine the number of homomorphisms to epimorphisms given in  \cref{Jones}, we simply subtract off the number of epimorphisms onto the subgroups of $G$ which we have determined already earlier in our list.

For step (3), for many triangle groups, there are only a few containments, so it is often fairly straightforward to reason when different groups have to be conjugate. We illustrate with an example.

\begin{exa}
\label{ex-c5}
The cyclic group $C_5$ acts with signature $(0;5,5,5)$ on a surface of genus $2$. Using \cref{Jones}, the number of homomorphisms from $\Gamma$ with signature $(0;5,5,5)$ onto $C_5$ is $12$. Since there are no non-trivial normal subgroups of $C_5$, all of these must be epimorphisms, and thus there are $12/|\Aut(C_5)|=12/4=3$ distinct surface subgroups of $\Gamma$ with quotient $C_5$. From Breuer's lists, we know there is no group action on a surface of genus $2$ with signature $(0;3,3,5)$. In particular, none of these surface subgroups are normal in the Fuchsian overgroup $\Gamma_1$ with signature $(0;3,3,5)$. Since $\Gamma$ is normal in $\Gamma_1$, the three surface subgroups must all be conjugate within $\Gamma_1$ (since $[\Gamma_1 :\Gamma]=3$). In particular, there is exactly one surface of genus $2$ up to conformal equivalence on which the group $C_5$ acts with signature $(0;5,5,5)$.
\end{exa}

We use this information to determine full automorphism groups in the following way.  Suppose we know a particular containment of Fuchsian triangle groups $\Gamma \leq \Gamma_1$, both of which appear in a list such as Breuer's with corresponding finite groups $G$ and $G_1$. Then if we know that $G_1$ contains a subgroup isomorphic to $G$, it is possible that $G$ is not the full automorphism group of a given surface on which $G$ acts. In order to decide this, we can count the number of distinct surfaces whose surface subgroup is contained in $\Gamma$ normally with quotient $G$, and compare that to the number of distinct surfaces whose surface subgroup is contained in $\Gamma_{1}$ normally with quotient $G_1$ -- if there are more subgroups inside $\Gamma$, then we know for sure there exists at least one surface $\X$ for which $G$ acts but does not extend to $G_1$. Of course, this doesn't mean it doesn't extend to other actions given by other Fuchsian group containments, but we can check each containment individually (and many are unique anyway).

Things are trickier if there are the same number or less than the same number of normal surface subgroups of $\Gamma$ with quotient $G$ than normal surface subgroups of $\Gamma_1$ with quotient $G_1$. The difficulty lies in the fact that it is not always guaranteed if $G_1$ acts on $\X$ then the subgroup $G$ corresponds to $\Gamma$ -- it could correspond to a completely different Fuchsian group which may not even be a triangle group. Fortunately, it is easy to determine signatures of Fuchsian subgroups of a given index using the following result of Singerman:

\begin{thm}
\label{Singerman1}Let $\Gamma_1$ have signature $(g_{\Gamma_1};m_{1},\dots ,m_{r})$. Then $\Gamma_1$ contains a subgroup $\Gamma$ of finite index with signature $(g_{\Gamma};n_{1,1},n_{1,2},\dots ,n_{1,\theta_{1}},\dots n_{r,\theta_{r}})$ if and only if

\begin{enumerate}

\item
There exists a finite permutation group $G$ transitive on $[\Gamma_1 :\Gamma]$ points and an epimorphism $\Phi \colon \Gamma_1 \rightarrow G$ such that the permutation $\Phi (\zeta_{j})$ has precisely $\theta_{j}$ cycles of length less then $m_{j}$, the lengths of these cycles being $$m_{j} /n_{j,1},\dots m_{j}/n_{j,\theta_{j}}.$$

\item
$[\Gamma_1 :\Gamma]=(2g_{\Gamma}-2+\sum\limits_{j=1}^{r} \sum\limits_{i=1}^{\theta_{j}} (1-\frac{1}{n_{i,\theta_{j}}} ))/(2g_{\Gamma_1}-2+\sum\limits_{i=1}^{r} (1-\frac{1}{m_{i}} ) )$.

\end{enumerate}

\end{thm}
Now, if $\Gamma$ is the unique subgroup of a given index, then each surface subgroup of $\Gamma_1$ with quotient $G_1$ will also be a surface subgroup of $\Gamma$ with quotient $G$. In particular, if there are the same number of surface subgroups, then all such surfaces have $G_1$ as their full automorphism group and not $G$.  If $\Gamma$ is not unique the problem becomes much more difficult, and we need to start looking at generating vectors for $G_1$ corresponding to each of the different surfaces to determine the signature of the subgroups corresponding to $G$ using  \cref{Singerman1}. Rather than present complete details, we illustrate with an explicit example below, and then use these techniques to explain how to classify full automorphism groups of superelliptic surfaces in the next section.

\begin{exa}
In  \cref{ex-c5}, we saw that there is exactly one surface of genus $2$ up to conformal equivalence on which the group $C_5$ acts with signature $(0;5,5,5)$. In Breuer's list for genus $2$, there is also an action by the cyclic group $C_{10}$ with signature $(0;2,5,10)$. Using \cref{Jones}, the number of homomorphisms from $\Gamma_1$ with signature $(0;2,5,10)$ onto $C_{10}$ is $4$. Since there are no non-trivial normal subgroups of $C_{10}$ containing an element of order $10$, all of these must be epimorphisms, and thus there is $4/|\Aut(C_{10})|=4/4=1$ distinct surface subgroup of $\Gamma_1$ with quotient $C_{10}$. Now, using \cref{Singerman1}, it is easy to see that $\Gamma_1$ has a unique subgroup of index $2$, this subgroup having signature $(0;5,5,5)$. In particular, it follows that corresponding unique surface of genus $2$ on which $C_5$ acts with signature $(0;5,5,5)$ has a larger automorphism group -- at least $C_{10}$. In particular, since this is the unique surface on which $C_5$ acts with this signature, it cannot possibly be the full automorphism group of a genus $2$ surface.

\end{exa}

\subsection{Superelliptic curves revisited}\label{subsec-revisit}

In order to illustrate the methods we have outlined for finding full automorphism groups, we conclude our discussion on compact Riemann surfaces by describing how to use the techniques we have described to find the full automorphism group of all superelliptic surfaces of level $n$ (henceforth just superelliptic surfaces). By {\em superelliptic surface}, we mean a compact Riemann surface $\X$ with a cyclic group of automorphisms $H$ of order $n$, called a {\em superelliptic group}, with the property that $\X/H$ has genus $0$ and that every branch point of the map $\pi_{H} \colon \X \rightarrow \X/H$ has order $n$. Such surfaces are natural generalizations of hyperelliptic surfaces, where $n=2$ whose full automorphism groups were completely classified in \cite{bgg}. It can be shown that, assuming the above definition of superelliptic, the curve has an equation of the form $y^n=f(x)$ with some restrictions on the factors of $f(x)$. See the references at the beginning  \cref{sec-super}.

\subsection{Preliminaries on Superelliptic Surfaces} Before we start our analysis, we shall first introduce some basic terminology, notation and facts about superelliptic surfaces from the point of view of Fuchsian groups. Henceforth, for our analysis, $\X$ will denote a superelliptic surface and $H$ a superelliptic group with generator $\tau$. Now since $\X /H$ has genus $0$ and every branch point of the map $\pi_{H} \colon \X \rightarrow \X/H$ has order $n$, it follows that $H$ has signature $(0;n,\dots ,n)$.

\begin{figure}[h]
        \centering
  $
\xymatrix @R=.4in @C=.4in {
\Gamma_{\Aut(\X)} \ar[r]_{\rho_{\Gamma_{\Aut(\X)}}} \ar[d]_{\bigtriangledown} & \Aut(\X) \ar[r]_{\rho_{K}} \ar[d]_{\bigtriangledown}  & K  \\
\Gamma_{H} \ar[r]_{\rho_{\Gamma_{\Aut(\X)}}} \ar[d]_{\bigtriangledown} & H \ar[r]_{\rho_{K}} \ar[d]_{\bigtriangledown} & 1 \\
\Lambda \ar[r]_{\rho_{\Gamma_{\Aut(\X)}}} & 1 &  \\
}
$
\caption{\label{Groups}Groups and Quotients}
 \end{figure}
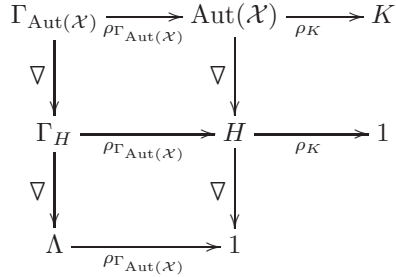

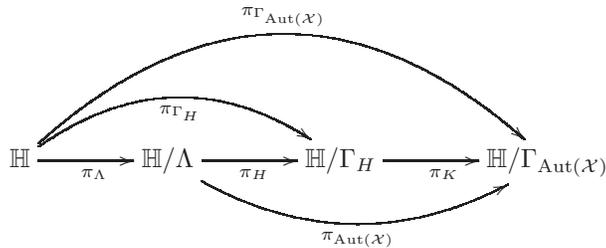
\begin{figure}[h]
        \centering
$
\xymatrix @R=.2in @C=.5in {
{\mathbb H}  \ar@/^4pc/[rrr]^{\pi_{\Gamma_{\Aut(\X)}}} \ar@/^2pc/[rr]_{\pi_{\Gamma_{H}}} \ar[r]_{\pi_{\Lambda}} & \ar@/_2pc/[rr]_{\pi_{\Aut(\X)}}  {\mathbb H} /\Lambda \ar[r]_{\pi_{H}} & {\mathbb H} /\Gamma_{H}  \ar[r]_{\pi_{K}} &  {\mathbb H} /\Gamma_{\Aut(\X)} \\
}
$
\caption{\label{Covers}Maps and Quotient Spaces}
  \end{figure}

Now, if $\Aut(\X) >H$, then using the notation we have introduced, after appropriate identifications, we have the tower of groups and epimorphisms illustrated in \cref{Groups} and corresponding to this, the tower of surfaces and holomorphic maps between them illustrated in \cref{Covers}.

Observe that the map $\pi_{K}$ is a finite Galois map with covering group $K$ from the Riemann sphere to itself, and all such maps are well known. We summarize the possibilities for $K$ and the ramification data of the map $\pi_{K}$ in  \cref{Standard}. The importance of these observations is that  we can use it to help find $\Aut(\X)$ and the possible signatures of $\Gamma_{\Aut(\X)}$. Specifically, since $H$ is normal in $\Aut(\X)$ with quotient $K$, $\Aut(\X)$ will be a group satisfying the short exact sequence $$1\rightarrow H \rightarrow \Aut(\X)\rightarrow K\rightarrow 1$$ and all such groups are relatively easy to find. For our purposes, we can take the more general list of such groups in any characteristic given in \cref{th14} and refine this list to the admissible $K$ for characteristic $p=0$. 

\begin{table}[h]
$
\begin{array}{||c|c||}
\hline \hline
{\rm{Group}} & {\rm{Branching \hspace{.1cm} Data}} \\
\hline \hline
C_{n} & (n,n) \\
\hline
D_{n} & (2,2,n) \\
\hline
A_{4} & (2,3,3) \\
\hline
S_{4} & (2,3,4) \\
\hline
A_{5} & (2,3,5) \\
\hline \hline
\end{array}
$
\caption{\label{Standard}Groups of Automorphisms of the Riemann Sphere and Branching Data}
\end{table}

\vspace{-5mm}

Now the signature of $\Gamma_{\Aut(\X)}$ depends only on which branch points of the map $\pi_{H}$ are also ramification points of the map $\pi_{K}$. We summarize (see for example  Proposition 3 in \cite{Woo1} for details).

\begin{thm}
\label{RamData}
The signature of $\Gamma_{\Aut(\X)}$ takes one of the following forms:

\begin{enumerate}

\item
If $K= C_{m}$ and $b_1,b_2\in \X /\Aut(\X)$ are the branch points of $\pi_K$, the signature of $\Gamma_{\Aut(\X)}$ is $(0;a_{1}m ,a_{2}m ,n,\dots ,n)$ where $a_i=n$ if $\pi_{K}^{-1}(b_i)$ contains a branch point of the map $\pi_{H}$ and $a_i=1$ else.

\item
If $K\neq C_{m}$ and $(d_{1},d_{2},d_{3})$ is the branching data of the quotient map $\pi_{K}$ with corresponding branch points $b_1,b_2,b_3\in \X /\Aut(\X)$ respectively, the signature of $\Gamma_{\Aut(\X)}$ is $(0;a_{1}d_{1} ,a_{2}d_{2} ,a_{3}d_{3} ,n,\dots ,n)$ where $a_i=n$ if $\pi_{K}^{-1}(b_i)$ contains a branch point of the map $\pi_{H}$ and $a_i=1$ else.

\end{enumerate}

\end{thm}

\subsection{Finding Full Automorphism Groups of Superelliptic Surfaces}

With the necessary terminology and notation introduced, we can now formalize the problem and explain how to solve it. First, we recall the problem: for a valid group signature pair ($G$,$\mathcal{S}$) for a superelliptic surface, we want to know if there exists a  superelliptic surface $\X$ on which $G$ acts with signature $\mathcal{S}$ as the full automorphism group of $\X$, or conversely, if for every  such superelliptic surface, $G$ is always contained in some larger group.

Two immediate observations. First, clearly if the signature $\mathcal{S}$ does not appear in Singerman's list, then there always exists a superelliptic surface $\X$ on which $G$ acts with signature $\mathcal{S}$ as the full automorphism group of $\X$. Therefore, to solve this problem, we just need to consider possible group signature pairs for which the signature appears in Singerman's list. Second, since we are assuming $H<G$, the orbit genus of the signature of $\Gamma_G$ must be $0$, so we can eliminate cases N1 through N4 from Singerman's list as possibilities for $\mathcal{S}$.

Next, we can eliminate further possibilities using the fact that our surface is superelliptic. Specifically, if an extension of $G'$ of $G$ exists then $H$ would also be central in $G'$ and would therefore be isomorphic to one of the groups appearing in \cref{th14} (with appropriate $K$) with corresponding quotient $K'=G'/H$ being an extension of $K=G/H$ with corresponding ramification data from \cref{Standard}. In particular, when $K=S_4$ or $A_5$, there are no possible extensions, so in either of these cases, $G$ acts as the full automorphism group. There are also restrictions to which groups a given $K$ can extend to for the three remaining cases, $A_4$, $D_m$ and $C_m$ which we shall use in our analysis below.

At this point, we have enough information to tabulate all possibilities for the signatures which might arise as signatures for non-maximal groups of automorphisms of superelliptic surfaces. We shall tabulate these signatures organizing them by isomorphism classes of $K$. For a given signature of $\Gamma_G$ in each such table, we also include the signature for any Fuchsian group $\Delta$ with $\Delta \geq \Gamma_{G}$ with signature coming from Singerman's list together with $K'$, the corresponding extension of $K$ which is needed to test for the maximality of $G$. The arguments for why we can eliminate each group and signature is a little different for each of the possible $K's$, so we break up our arguments accordingly.

Suppose first that $K=C_{m}$ and suppose that $K'$ is some extension corresponding to $G'$, an extension of $G$. Now, it is not possible for $K=C_{ms}$ since using \cref{Singerman1}, the number of branch points of $\Gamma_{G}$ would be strictly bigger than that of $\Gamma_{G'}$, and so the signatures would not appear in Singerman's list. Next, it is always possible for $K'=D_{m}$ with corresponding signature $(0;2,2,m)$. Moreover, if $K$ is contained in some larger dihedral group, it will necessarily also be contained in $D_{m}$, so to check for maximality, we just need to check extension to $D_{m}$. For $K'=S_4$ or $A_5$, each cyclic subgroup $C_{m}$ is also a subgroup of $D_{m}$ within $K'$, so in particular, to check for maximality, we just need to check extension to $D_{m}$. Finally, for $K=A_4$, any $C_2$ is contained in $D_{2}$, so to check for maximality, we just need to check extension to $D_{2}$. For $C_{3}\leq A_4$ however, there are no intermediate subgroups, so to check for maximality, we need to check extension to $A_4$.

\begin{table}[h]
\begin{tabular}{||c|c|c|c|c||}
\hline \hline
$K$ & Signature of $\Gamma_{G}$ & $K'$ & Signature of $\Delta$ & Further Conditions \\
\hline \hline
$C_{m}$ & $(m,m,n,n)$ &$D_{m}$ & $(0;2,2,m,n)$ & \\
\hline
$C_{m}$ & $(mn,mn,n,n)$ & $D_{m}$ & $(0;2,2,mn,n)$ & \\
\hline
$C_{m}$ & $(m,m,n)$ & $D_{m}$ & $(0,2,m,2n)$ & \\
\hline
$C_{3}$ & $(3,n,3n)$ & $A_{4}$ & $(0;2,3,3n)$ & $m=3$ \\
\hline
$C_{m}$ & $(mn,mn,n)$ &  $D_{m}$ & $(0;2,mn,2n)$ & \\
\hline \hline
\end{tabular}
\caption{\label{nonmaxsigscyclic}Potential Non-Maximal Signatures with $K=C_{m}$}
\end{table}

\vspace{-5mm}

Now suppose that $K=D_{m}$ and suppose that $K'$ is some extension corresponding to $G'$, an extension of $G$. Using \cref{Singerman1} to find the corresponding signature of the surface group for $\Lambda$ in $\Gamma_{G}$ and $\Gamma_{G'}$, we see that if $K'=D_{ms}$, then either $s=2$ or $s=4$. In particular, in each case, it will necessarily be contained in $D_{2m}$, so to check for maximality, we just need to check extension to $D_{2m}$. The only dihedral subgroup of $K'=A_4$ is $D_2$ with $m=2$,  and in this case we need to check for extension to $A_4$. For $K'=S_4$, there are three different dihedral subgroups -- $D_2$, $D_3$ and $D_4$. If a given $D_2$ extends to $S_4$, then it also extends to $D_4$, so to check maximality, we only need to check maximality in $D_4$. The group $D_3$ has index $4$ in $S_4$. Neither of the signatures $T9$ or $T10$ with index $4$ are valid signatures for $D_3$ and $S_4$ actions, so in this case we do not need to check maximality. Finally, the group $D_4$ has index $3$ in $S_4$. Of the two signatures in Singerman's list with index $3$, only $T11$ is valid for $D_4$ and $S_4$ with $(0;2,2n,4n)$ and $(0;2,3,4n)$ respectively, so this is the only case where we need to check maximality of $D_4$ in $S_4$.

\begin{table}[h]
\begin{tabular}{||c|c|c|c|c||}
\hline \hline
$K$ & Signature of $\Gamma_{G}$ & $K'$ & Signature of $\Delta$ & Further Conditions \\
\hline \hline
$D_{m}$ & $(0;2n,2n,n,n)$ & $D_{2m}$ & $(0;2,2,n,2n)$ & $m=n$\\
\hline
$D_{m}$ & $(0;2n,2n,m)$ &$D_{2m}$ & $(0;2,2n,2m)$ & \\
\hline
$D_{m}$ & $(0;2n,2n,mn)$ &$D_{2m}$ & $(0;2,2n,2nm)$ & \\
\hline
$D_{2}$ & $(0;2n,2n,2n)$ & $A_4$ & $(0;2n,3,3)$ & $m=2$ \\
\hline
$D_{4}$ & $(0;2,2n,4n)$ & $S_4$ & $(0;2,3,4n)$ & $m=4$ \\
\hline \hline
\end{tabular}
\caption{\label{nonmaxsigsdihedral}Potential Non-Maximal Signatures}
\end{table}

\vspace{-5mm}

Lastly, now suppose that $K=A_4$ and suppose that $K'$ is some extension corresponding to $G'$, an extension of $G$. The only possibilities for $K'$ are $S_4$ and $A_5$, so we look at each of these cases individually. If $K'=S_4$, then $A_4$ has index $2$ in $S_4$, and in this case there are multiple signatures from Singerman's list of the form $N8$ which we need to check to determine maximality. When $K'=A_5$, then $K=A_4$ is non-normal and has index $6$ in $K'$. This leaves $T7$ and $T8$ as the only possible signature pairs, but neither of these work for the inclusion of $A_4$ in $A_5$. We summarize in \cref{nonmaxsigsa4}.

\begin{table}[h]
\begin{tabular}{||c|c|c|c|c||}
\hline \hline
$K$ & Signature of $\Gamma_{G}$ & $K'$ & Signature of $\Delta$ & Further Conditions \\
\hline \hline
$A_{4}$ & $(0;2,3n,3n)$ & $S_4$ & $(0;2,3n,4)$ & \\
\hline
$A_{4}$ & $(0;2n,3,3)$ & $S_4$ & $(0;2,3,4n)$ & \\
\hline
$A_{4}$ & $(0;2n,3n,3n)$ & $S_4$ & $(0;2,3n,4n)$ & \\
\hline \hline
\end{tabular}
\caption{\label{nonmaxsigsa4}Potential Non-Maximal Signatures}
\end{table}

\vspace{-5mm}

At this point, we now have enough information that to finish the problem we can either apply \cref{Cond1} and \cref{Cond2} to the different possible generating vectors for each group and signature pair, or alternatively apply the counting methods developed in \cref{count}. Rather than go through every individual case, we illustrate with an explicit example. The other signature pairs given in \cref{nonmaxsigscyclic}, \cref{nonmaxsigsdihedral} and \cref{nonmaxsigsa4} yield results with similar congruence conditions given in \cref{ex-cyclic}.

\begin{exa}\label{ex-cyclic}
Consider the signature $(0;mn,mn,n)$. Since the signature has periods of order $mn$, the only possible quotient group in this case is $C_{nm}$. If $C_{nm}=\langle x\rangle$, then after appropriate automorphism of $C_{nm}$, any $(0;mn,mn,n)$-generating vector of $C_{nm}$ will have the form $(x,x^a,x^{nm-a-1})$ where $(a,nm)=1$ and $(nm,nm-a-1)=m$.  

Applying $(2)$ of \cref{Cond1}, this group extends to a group with signature $(0; 2, mn$, $2n)$ if and only if the maps $x\rightarrow x^a$ and $x^a\rightarrow x$ induce an automorphism of $C_{nm}$ (note that the condition on the third generator is trivially satisfied since the group is Abelian). However, this happens only if $x^{a^2}=x$ i.e. $a^2\equiv 1\mod nm$. In particular, unless $a^2\equiv 1\mod nm$ for every $a$ satisfying the two different congruences, then there exists at least one surface on which $C_{nm}$ acts as the full automorphism group with signature $(0;mn,mn,n)$
\end{exa}

We note that there are examples from \cref{ex-cyclic} for which the group is never maximal, such as when $n=4$ and $m=10$, and when there exist maximal actions such as $m=n=8$.

\nocite{*}

\bibliographystyle{amsalpha}
\bibliography{ref}{}

\end{document}